\documentclass[a4paper,french]{rnti}

\usepackage{dsfont}

\usepackage{graphicx}

\usepackage{amsfonts,amsmath, epsfig,amssymb,array}
\newcommand{\bx}{\mathbf{x}}

\newcommand{\bz}{\mathbf{z}}

\newcommand{\bt}{\mathbf{t}}
\newcommand{\bT}{\mathbf{T}}
\newcommand{\bW}{\mathbf{W}}
\newcommand{\bsPhi}{\boldsymbol{\Phi}}
\newcommand{\bsw}{\boldsymbol{w}}

\newcommand{\bsv}{\boldsymbol{v}}
\newcommand{\bsbeta}{\boldsymbol{\beta}}
\newcommand{\btheta}{\boldsymbol{\theta}}

\titrecourt{Modèle à processus latent et algorithme EM pour la régression non linéaire}

\nomcourt{Chamroukhi et al.}

\titre{Modèle à processus latent et algorithme EM pour la régression non linéaire}
\auteur{Faicel Chamroukhi \affil{1}\affilsep\affil{2}, \ Allou Samé\affil{1}, \ Gérard Govaert\affil{2}, \ Patrice Aknin\affil{1} }
\affiliation{
    \affil{1}Institut National de Recherche sur les Transports et leur Sécurité\\
    Laboratoire des Technologies Nouvelles\\
 2 rue de la Butte Verte, 93166 Noisy-le-Grand Cedex\\
          \{chamroukhi, same, aknin\}@inrets.fr \\ \\
    \affil{2}Université de Technologie de Compiègne\\
     Laboratoire HEUDIASYC, UMR CNRS 6599 \\
BP 20529, 60205 Compiègne Cedex\\
          govaert@utc.fr\\
 }

\resume{Cet article propose une méthode de régression non linéaire
qui s'appuie sur un modèle intégrant un processus latent qui permet
d'activer préférentiellement, et de manière souple, des sous-modèles
de régression polynomiaux. Les paramètres du modèle sont estimés par
la méthode du maximum de vraisemblance mise en {\oe}uvre par un
algorithme EM dédié. Une étude expérimentale menée sur des données
simulées et des données réelles met en évidence de bonnes
performances de l'approche proposée.}

\summary{ A non linear regression approach which consists of a
specific regression model incorporating a latent process, allowing
various polynomial regression models to be activated preferentially
and smoothly, is introduced in this paper. The model parameters are
estimated by maximum likelihood performed via a dedicated EM
algorithm. An experimental study using simulated and real data sets
reveals good performances of the proposed approach.}

\begin{document}

\section{Introduction}

La régression non linéaire est un problème central dans de nombreux
domaines qui concer\-nent la prédiction, le débruitage de signaux et
leur paramétrisation. Son but est de caractériser au mieux la
relation (non linéaire) existant entre une variable dépendante
(grandeur physique), que nous supposerons scalaire dans cet article,
et une variable indépendante qui, comme c'est souvent le cas, est
liée au temps. Cette relation non linéaire peut être due au fait que
les données sont issues d'un modèle physique intrinsèquement non
linéaire par rapport au temps, ou que le processus de génération des
données comporte différents régimes linéaires (voire polynomiaux)
qui se succèdent au cours du temps.

Plusieurs modèles ont déjà été proposés dans le cadre de
l'apprentissage statistique pour résoudre ce type de problème. Parmi
ces approches, on peut citer les modèles polynomiaux par morceaux
\citep{McGee}, \citep{bellman, stone}, \citep{sameSFC2007}, les méthodes à base de B-splines \citep{deboor1978, bishopPRML}, le perceptron multi-couche dans sa version
régressive \citep{bishopPRML} et les méthodes à fonctions de base
radiale \citep{bishopPRML}. La plupart de ces méthodes ramènent
généralement le problème de régression non linéaire à des problèmes
de régression linéaire simples à résoudre.

Dans ce travail, nous proposons une méthode alternative qui consiste
à remplacer le modèle de régression non linéaire habituel par un
modèle de régression intégrant un processus caché, qui
permet d'activer préférentiellement et de manière souple un modèle de régression
polynomial parmi $K$ modèles. L'utilisation d'une fonction
logistique comme loi conditionnelle des variables latentes assure
une souplesse de transition (lente ou rapide) entre les différents
polynômes, ce qui permet d'obtenir une modélisation correcte des non
linéarités. La loi conditionnelle de la variable dépendante, sous ce
modèle, constitue un mélange d'experts \citep{jacobsME},
\citep{jordanHME} contraint, avec une variance commune pour toutes
les composantes du mélange. Une estimation des paramètres du modèle
par l'algorithme EM est ainsi proposée.

Cette méthode est exploitée dans le cadre d'une application de suivi d'état de fonctionnement du mécanisme d'aiguillage des rails, où nous avons été amenés à paramétriser des signaux non linéaires acquis durant des man{\oe}uvres d'aiguillage.

Dans la section 2 nous décrivons comment le modèle de régression à
processus latent peut être utilisé dans le cadre de la régression
non linéaire. Ensuite nous présentons la méthode d'estimation des
paramètres via l'algorithme EM. La section 4 montre les performances
de la méthode proposée sur des données simulées et la section 5 est
consacrée à son application sur des données réelles issues du
domaine ferroviaire.

\section{Régression non linéaire et modèle à processus latent}

\subsection{Cadre général}

On suppose disposer d'un échantillon $((x_1,t_1),\ldots,(x_n,t_n))$,
où $x_i$ désigne la variable aléatoire scalaire dépendante et $t_i$
la variable temporelle indépendante. Le problème de régression non
linéaire, sous sa formulation générique \citep{reg_non_lin},
consiste à estimer une fonction $f$ de paramètre $\btheta$,  en
considérant le modèle
\begin{equation}
M1 \ :\   x_i = f(t_i;\btheta) + \varepsilon_i,
\label{eq. M1}
\end{equation}
où les $\varepsilon_i$ sont des bruits gaussiens centrés $i.i.d$ de
variance $\sigma^2$. L'estimation s'effectue généralement par la
maximisation de la vraisemblance ou de manière équivalente par la
minimisation du critère des moindres-carrés donné par
\begin{equation}
C_1(\btheta)=\sum_{i=1}^n (x_i - f(t_i;\btheta))^2.
\label{eq. C1}
\end{equation}
Dans le modèle $M1$, la fonction $f(t;\btheta)$ représente
l'espérance de $x$ conditionnellement à $t$. La minimisation du
critère $C_1$, sous certaines conditions de régularité sur $f$
(\citet{reg_non_lin}), fournit un estimateur asymptotiquement
gaussien, efficace et sans biais du paramètre $\btheta$. En
pratique, on a souvent recours à
des algorithmes itératifs d'optimisation qui convergent localement.

Pour couvrir un panel assez large de fonctions  non linéaires de
régression qui soient facilement paramétrables, nous optons pour des
fonctions qui peuvent s'écrire sous la forme d'une somme finie de
polynômes pondérés par des fonctions logistiques :
\begin{equation}
f(t;\btheta) = \sum_{k=1}^K \pi_k(t;\bsw) \bsbeta_k^T\bt,
\label{modele de la fonc non lineaire}
\end{equation}
avec
\begin{equation}
\pi_k(t;\bsw)=\frac{\exp(w_{k 0}+w_{k 1}
t)}{\sum_{\ell=1}^K \exp(w_{\ell 0}+w_{\ell 1} t)},
\label{softmax}
\end{equation}
où le vecteur $\bsbeta_k=(\beta_{k0},\ldots,\beta_{kp})^T$ de
$\mathds{R}^{p+1}$ désigne l'ensemble des coefficients d'un polynôme de
degré $p$ et $\bt=(1,t,t^2,\ldots,t^p)^T$ son vecteur de monômes
associés. Le vecteur $\bsw=(w_{10},w_{11},\ldots,w_{K0},w_{K1})$ de
$\mathds{R}^{2K}$ désigne l'ensemble des paramètres de la fonction
logistique et $\btheta=(\bsw,\bsbeta_1,\ldots,\bsbeta_K)$ représente
l'ensemble des paramètres de la fonction $f$. 

 Ce modèle particulier de régression peut s'interpréter comme étant
un modèle formé de différents sous-modèles de régression activés de
manière souple ou brusque par des fonctions logistiques.  Cette
méthode permet aussi de détecter les points de rupture d'un signal par le suivi au cours du
temps du processus latent, et de mesurer la vitesse d'évolution des proportions du mélange (souple ou brusque)
gr\^ace à la flexibilité de la transformation logistique utilisée. L'illustration de cette flexibilité se trouve dans \citep{chamroukhi_et_al_ESANN2009}.

Pour le modèle $M1$, défini par les équations (\ref{eq. M1}), (\ref{modele de la fonc non lineaire}) et (\ref{softmax}), le choix particulier de la fonction $f$ ne permet pas d'obtenir une solution analytique du problème de minimisation du critère $C_1$. Il faut s'appuyer sur une procédure numérique d'optimisation du type Gauss-Newton, Newton-Raphson ou quasi-Newton qui converge localement.

\subsection{Modèle de régression à processus latent proposé}

Pour atteindre l'objectif d'estimation de la fonction de
régression $f$ du modèle $M1$, nous proposons d'utiliser une méthode
alternative s'appuyant sur un modèle génératif intégrant un
processus latent discret $(z_1,\ldots,z_n)$ avec $z_i\in
\{1,\ldots,K \}$. Ce modèle est défini par :
\begin{equation}
M2 \ : \  x_i = \sum_{k=1}^K z_{ik}\bsbeta^T_k \bt_i + \varepsilon_i
\quad ; \quad \varepsilon_{i} \sim \mathcal{N} (0,\sigma^2),
\label{M2}
\end{equation}
où $z_{ik}$ vaut 1 si $z_i=k$ et vaut 0 sinon, et où les $\varepsilon_i$ sont
supposés être indépendants. Les variables $z_i$ du processus latent,
conditionnellement aux instants $t_i$, sont supposées être générées
indépendamment suivant la loi multinomiale
$\mathcal{M}(1,\pi_1(t_i;\bsw),\ldots,\pi_K(t_i;\bsw))$, où
$\pi_k(t_i;\bsw)$ est défini par l'equation (\ref{softmax}). \`A partir du modèle $M2$, on peut vérifier que conditionnellement à un sous-modèle de régression $k$ et au temps $t_i$, $x_i$ est distribué suivant une loi normale de moyenne $\bsbeta_k^T \bt_i$ et de variance $\sigma^2$
\begin{equation}
p(x_i|z_i=k,t_i;\bsPhi) = \mathcal{N}(x_i;\bsbeta_k^T \bt_i,\sigma^2),
\end{equation}
où $\bsPhi=(\btheta,\sigma^2)$ et $\mathcal{N}(\cdot;\mu,\sigma^2)$
désigne la fonction de densité d'une loi normale d'espérance $\mu$
et de variance $\sigma^2$. On peut montrer alors que conditionnellement à $t_i$, la variable aléatoire
$x_i$ est distribuée suivant le mélange de densités normales
\begin{eqnarray}
p(x_{i}|t_i;\bsPhi) &=
&\sum_{k=1}^K\pi_{k}(t_i;\bsw)\mathcal{N}\big(x_{i};\bsbeta^T_k\bt_{i},\sigma^2\big).
\label{melange}
\end{eqnarray}
 Il faut noter que cette dernière loi
conditionnelle est un mélange d'experts \citep{jacobsME},
\citep{jordanHME} contraint, avec une variance commune pour toutes
les composantes du mélange. Le modèle $M2$ constitue un modèle parcimonieux du modèle plus général supposant des variances différentes. En pratique, pour les données traitées,  ces deux modèles ont conduit à des ajustements quasi identiques.

On peut alors vérifier que le modèle de régression proposé conduit à
la même espérance conditionnelle que celle du modèle $M1$:
\begin{eqnarray}
E[x|t;\bsPhi] &=& \int_{\mathds{R}}x{}p(x|t;\bsPhi)dx  \nonumber \\
&=& \sum_{k=1}^K\pi_{k}(t;\bsw) \int_{\mathds{R}}x \mathcal{N}\big(x ; \bsbeta^T_k\bt,\sigma^2\big) dx\nonumber\\
&= &\sum_{k=1}^K \pi_k(t;\bsw) \bsbeta_k^T\bt \nonumber \\
&=& f(t;\btheta).
\label{eq. somme pond de polynomes}
\end{eqnarray}
Ainsi, grâce aux propriétés asymptotiques classiques de normalité,
d'absence de biais et d'efficacité de l'estimateur du maximum de
vraisemblance de $\btheta$, la fonction de régression $f$ estimée à
partir du modèle $M2$ est asymptotiquement identique à celle obtenue
à partir du modèle $M1$. Ce qui nous conforte sur le fait que le
modèle proposé peut être une bonne alternative pour résoudre le
problème de régression non linéaire, si on dispose d'un algorithme
adapté pour l'estimation de ses paramètres.

La section suivante montre comment les paramètres du modèle peuvent
être estimés par la méthode du maximum de vraisemblance.

\section{ Estimation des paramètres via l'algorithme EM}

Dans cette section, compte tenu du fait que la densité de la loi
conditionnelle de $x$ s'écrit sous la forme d'un mélange de
densités, nous exploitons le cadre élégant de l'algorithme EM
\citep{dlr} pour estimer ses paramètres.

Les hypothèses d'indépendance des $\varepsilon_i$ et d'indépendance
des $z_i$ conditionnellement aux $t_i$ entraînent l'indépendance des
$x_i$ conditionnellement aux $t_i$. La log-vraisemblance à maximiser
s'écrit donc
\begin{eqnarray}
L(\bsPhi) &=& \sum_{i=1}^{n}\log p(x_i|t_i;\bsPhi) \nonumber \\
 &=&\sum_{i=1}^{n}\log\sum_{k=1}^K\pi_{k}(t_i;\bsw)\mathcal{N}\big(x_{i};\bsbeta^T_k\bt_i,\sigma^2\big)\cdot
\label{vraisemblance}
\end{eqnarray}
Cette maximisation ne pouvant pas être effectuée analytiquement,
nous nous appuyons sur l'algorithme EM \citep{dlr} pour l'effectuer.
L'algorithme EM, dans cette situation, itère à partir d'un paramètre
initial  $\bsPhi^{(0)}$ les deux étapes suivantes jusqu'à la
convergence.


\subsection{Étape E (Espérance)} Cette étape consiste à calculer
l'espérance de la log-vraisemblan\-ce complétée \begin{equation}\log
p(\bx,\bz;\bsPhi)=\sum_{i=1}^{n}\sum_{k=1}^K  z_{ik} \log
[\pi_{k}(t_i;\bsw)\mathcal{N}
(x_{i};\bsbeta^T_k\bt_{i},\sigma^2)],\end{equation}
conditionnellement aux données observées et au paramètre courant
$\bsPhi^{(q)}$ ($q$ \'etant l'itération courante). Dans notre
situation, cette espérance conditionnelle s'écrit:
\begin{eqnarray}
Q(\bsPhi,\bsPhi^{(q)}) &=& E\left[
\log p(\bx,\bz;\bsPhi)|\bx;\bsPhi^{(q)}\right]\nonumber \\
&=&  \sum_{i=1}^{n}\sum_{k=1}^K E (z_{ik}|x_i;\bsPhi^{(q)}) \log [\pi_{k}(t_i;\bsw)\mathcal{N}
(x_{i};\bsbeta^T_k\bt_{i},\sigma^2)] \nonumber \\
&=&  \sum_{i=1}^{n}\sum_{k=1}^K p (z_{ik}=1|x_i;\bsPhi^{(q)}) \log [\pi_{k}(t_i;\bsw)\mathcal{N}
(x_{i};\bsbeta^T_k\bt_{i},\sigma^2)] \nonumber \\
&=&\sum_{i=1}^{n}\sum_{k=1}^K\tau^{(q)}_{ik}\log [\pi_{k}(t_i;\bsw)\mathcal{N}
(x_{i};\bsbeta^T_k\bt_{i},\sigma^2)],
\label{eq. fonction Q}
\end{eqnarray}
où
\begin{eqnarray}
 \tau^{(q)}_{ik} &=& p(z_{ik}=1|x_i;\bsPhi^{(q)}) \nonumber \\
 &=&\frac{\pi_{k}(t_i;\bsw^{(q)})\mathcal{N}(x_{i};\bsbeta^{T(q)}_k\bt_{i},\sigma^{2(q)})}
{\sum_{\ell=1}^K\pi_{\ell}(t_i;\bsw^{(q)})\mathcal{N}(x_{i};\bsbeta^{T(q)}_{\ell}\bt_{i},\sigma^{2(q)})},
\label{eq.tauik}
\end{eqnarray}
est la probabilité a posteriori que $x_i$ soit issu de la $k$ième
composante du mélange. Cette étape nécessite simplement le calcul
des $\tau^{(q)}_{ik}$.

\subsection{Étape M (Maximisation)}

Cette étape (de mise à jour) consiste à calculer le paramètre $\bsPhi^{(q+1)}$ qui maximise
$Q(\bsPhi,\bsPhi^{(q)})$ par rapport à $\bsPhi$. La quantité $Q(\bsPhi,\bsPhi^{(q)})$ s'écrit sous la forme d'une somme des deux quantités suivantes :
\begin{equation}
Q_1(\bsw)=\sum_{i=1}^{n}\sum_{k=1}^K \tau^{(q)}_{ik}\log
\pi_{k}(t_i;\bsw)
\end{equation}
et
\begin{equation}
Q_2(\bsbeta_1,\ldots,\bsbeta_K,\sigma^2) =
\sum_{i=1}^{n}\sum_{k=1}^{K} \tau^{(q)}_{ik}\log
\mathcal{N} (x_{i};\bsbeta^T_k\bt_{i},\sigma^2).
\end{equation}
Pour maximiser $Q(\bsPhi,\bsPhi^{(q)})$, il suffit donc de maximiser séparément les quantités $Q_1(\bsw)$ et $Q_2(\bsbeta_1,\ldots,\bsbeta_K,\sigma^2)$.

\medskip
La maximisation de $Q_2(\bsbeta_1,\ldots,\bsbeta_K,\sigma^2)$ est celle classique qu'on rencontre dans un mélange de lois gaussiennes d'espérance décrite par un modèle linéaire. On obtient les $\bsbeta^{(q+1)}_k$ en résolvant analytiquement $K$
problèmes de moindres-carrés ordinaires pondérés par les
$\tau^{(q)}_{ik}$:
\begin{eqnarray}
{ \bsbeta }_k^{T(q+1)} &=& \arg \min \limits_{\substack {\bsbeta}} \sum_{i=1}^{n} \tau^{(q)}_{ik} ( x_i-\bsbeta^{T}\bt_i )^2 \nonumber \\
 &=& (\bT^T\bW_k^{(q)}\bT)^{-1}\bT^T\bW_k^{(q)}\bx,
\label{estimation beta}
\end{eqnarray}
où $\bT$ est la matrice de régression de dimension $[n \times
(p+1)]$ définie par $$\bT=\left[\begin{array}{ccccc}1&t_{1}&t_{1}^2&\ldots&t_{1}^p\\
1&t_{2}&t_{2}^2&\ldots&t_{2}^p \\
\vdots&\vdots&\vdots&\vdots&\vdots\\
1&t_{n}&t_{n}^2&\ldots&t_{n}^p\end{array}\right],$$ $\bW_k^{(q)}$
est une matrice diagonale de dimension $[n \times n]$ ayant pour
éléments diagonaux \linebreak
$(\tau_{1k}^{(q)},\ldots,\tau_{nk}^{(q)})$ et $\bx =
(x_1,\ldots,x_n)^T$ est le vecteur de dimension $[(n+1)\times1]$ des
observations. La variance $\sigma^{2(q+1)}$, qui est identique pour
toutes les composantes du mélange, est donnée par :
\begin{eqnarray}
\sigma^{2(q+1)} &=&\arg \min \limits_{\substack{\sigma^2}} \left[n \log{\sigma^{2}} + \frac{1}{\sigma^{2}}\sum_{i=1}^{n} \sum_{k=1}^{K}\tau^{(q)}_{ik} (x_i-\bsbeta_k^{T}\bt_i )^2\right] \nonumber \\
&=& \frac{1}{n}\sum_{i=1}^{n}\sum_{k=1}^K\tau^{(q)}_{ik} (x_i-{\bsbeta}_k^{T(q+1)}\bt_i)^2.
\label{estimation sigma}
\end{eqnarray}

\medskip
La maximisation de $Q_1$ par rapport à $\bsw$ est un problème
convexe de régression logistique multinomial pondéré par les
$\tau^{(q)}_{ik}$. Contrairement au cas précédent, cette maximisation ne peut pas s'effectuer de manière analytique. On a recours à un algorithme, adapté à ce type de problème, et qui est lui-même itératif : l'algorithme IRLS (Iterative Reweighted Least Squares) \citep{irls}. Le paragraphe suivant décrit cet algorithme.

\paragraph{ Algorithme IRLS (Iteratively Reweighted Least Squares) : }
\label{paragraph: IRLS}

l'algorithme IRLS est utilisé pour maximiser $Q_1(\bsw)$ par rapport à $\bsw$  à l'étape M de chaque itération $q$ de l'algorithme EM. Puisque $\sum_{k=1}^{K} \pi_k(t_i;\bsw)=1$, les valeurs de $w_{K0}$ et $w_{K1}$ sont fixées à zéro pour éviter les problèmes d'identification.
L'algorithme IRLS est équivalent à l'algorithme de Newton-Raphson
qui consiste à partir d'un vecteur paramètre initial $\bsw^{(0)}$ et
à s'appuyer, à chaque nouvelle itération $c+1$,  sur la formule
suivante :
\begin{equation}
\bsw^{(c+1)}=\bsw^{(c)}-\left[{H(\bsw^{(c)})}\right]^{-1}g(\bsw^{(c)}),
\label{eq.IRLS}
\end{equation}
où $ H(\bsw^{(c)})$ et $g(\bsw^{(c)})$ sont respectivement la hessienne et le gradient de $Q_1(\bsw)$ calculés avec le paramètre $\bsw^{(c)}$. La matrice hessienne $H(\bsw^{(c)})$ est formée de
$(K-1)\times(K-1)$ matrices blocs $H_{k\ell}(\bsw^{(c)})$ ($k,\ell=1,\ldots,K-1$)
\citep{chen99}, où
\begin{eqnarray}
 H_{k\ell}(\bsw^{(c)})=-\sum_{i=1}^{n} \pi_k(t_i;\bsw^{(c)}) [\delta_{k\ell} - \pi_{\ell}(t_i;\bsw^{(c)})] \bsv_i {\bsv_{i}}^T,
\end{eqnarray}
$\delta_{k\ell}$ étant le symbole de Kronecker ($\delta_{k\ell}$ = 1
si $k=\ell$, 0 sinon) et $\bsv_i=(1,t_i)^T$. Le gradient de $Q_1(\bsw)$ s'écrit :
\begin{eqnarray}
g(\bsw^{(c)}) &=& [g_{1}(\bsw^{(c)}),\ldots,g_{K-1}(\bsw^{(c)})]^{T},
\end{eqnarray}
avec
\begin{eqnarray}
g_{k}(\bsw^{(c)}) &=& \sum_{i=1}^{n} [\tau^{(q)}_{ik} - \pi_{k}(t_i;\bsw^{(c)})]
\bsv^{T}_i \quad \forall \quad k=1,\ldots,K-1.
\end{eqnarray}
L'équation (\ref{eq.IRLS}), initialisée avec le paramètre
$\bsw^{(q)}$, fournit à la convergence le paramètre $\bsw^{(q+1)}$.

Chaque itération de l'algorithme EM proposé fait cro\^itre le critère de log-vraisemblance \citep{mclachlan_and_krishnan_EM}. Cela résulte de la maximisation de $Q(\bsPhi,\bsPhi^{(q)})$ et de l'inégalité de Jensen \citep{jensen}. Il a été prouvé par \citep{jordan_and_xu_1995} et \citep{xu_and_jordan_1996}, notamment pour le cas du mélange gaussien et celui du mélange d'experts, que l'algorithme EM permet de maximiser localement la log-vraisemblance. On peut limiter le nombre d'itérations de la procédure IRLS interne à l'algorithme EM. Cette version consiste à faire cro\^itre $Q(\bsPhi,\bsPhi^{(q)})$ à chaque itération au lieu de la maximiser. On peut par exemple limiter le nombre d'itérations de l'IRLS  jusqu'à une seule itération \citep{foulley2002}. On obtient ainsi un algorithme EM généralisé (GEM) \citep{dlr, mclachlan_and_krishnan_EM}, qui possède les mêmes propriétés de convergence que l'algorithme EM. En pratique, nous avons pu observer que cette limitation entrainait une augmentation du nombre d'itérations de l'algorithme EM. Par conséquent, nous avons opté pour une stratégie qui consiste à initialiser aléatoirement la procédure IRLS seulement pour la première itération de l'algorithme EM. Pour cette première initialisation, la convergence de l'algorithme IRLS requiert une quinzaine d'itérations. \`A partir de la deuxième itération de l'algorithme EM, l'algorithme IRLS défini par l'équation (\ref{eq.IRLS}) est initialisé avec le paramètre $\bsw^{(q)}$ estimé à l'itération précédente de l'algorithme EM. Au delà de la quatrième itération de l'algorithme EM, on observe que l'algorithme IRLS converge en moins de 5 itérations. Les critères d'arr\^et utilisés pour l'algorithme EM et l'algorithme IRLS consistent en  des valeurs seuils sur les variations relatives des log-vraisemblances à maximiser ($ |\frac{L(\boldsymbol{\Phi}^{(q+1)})-L(\boldsymbol{\Phi}^{(q)})}{L(\boldsymbol{\Phi}^{(q)})}|<10^{-6} $ pour EM  et $|\frac{Q_1({\bsw}^{(c+1)})-Q_1({\bsw}^{(c)})}{Q_1({\bsw}^{(c)})}|<10^{-6}$ pour IRLS) ainsi qu'à des  nombres d'itérations maximum ($1000$ pour EM et $50$ pour IRLS).

 \medskip

\subsection{Choix du nombre de composantes et de l'ordre des polynômes de régression}

Les valeurs optimales du nombre de composantes $K$ du modèle et de
l'ordre $p$ des polynômes de régression peuvent être obtenues en
maximisant le critère d'information bayésien (BIC) \citep{BIC}
défini comme étant le critère de vraisemblance obtenu à la
convergence de l'algorithme, pénalisé par un terme qui dépend du
nombre de paramètres libres du modèle. La pénalisation utilisée ici est $pen = -\frac{\nu(K,p) \log(n)}{2}$ où $\nu(K,p) = K(p+3)-1$ est le nombre total de paramètres libres à estimer.

Il est possible d'utiliser aussi le critère d'entropie normalisée NEC (Entropy Normalized Criterion) \citep{celeux_and_soromenho_NEC, biernacki_et_al_improvement_NEC}  adapté au contexte des modèles de mélange.

 \medskip
\subsection{Région de confiance de la courbe de régression $f(t;\btheta) $}

\`A partir de la modélisation proposée, on peut montrer, comme dans le cadre des mélanges hiérarchiques d'experts pour les modèles linéaires généralisés,  que sous certaines conditions de régularité,  $f(t;\hat{\btheta})$ suit asymptotiquement une loi normale de moyenne   $f(t;\btheta)$ et de variance $s^2 (t;\bsPhi)$ \citep{Jiang_and_tanner_IEEEinft_99}. Cependant, comme les ``vrais'' paramètres  $\bsPhi$ et $\btheta$ sont inconnus, on propose de les remplacer par leurs estimations $\hat{\bsPhi}$ et $\hat{\btheta}$ fournies par l'algorithme EM. Dans le cas de la régression linéaire, la région de confiance approchée au niveau $1-\alpha$ est définie par :
\begin{equation}        
f(t;\hat{\btheta}) \pm  \sqrt{\nu_{\btheta} p(\nu_{\btheta},n-\nu_{\btheta}; \alpha)} s(t;\hat{\bsPhi})\;,
\label{eq. region de confiance cas lineaire}
 \end{equation}
où $\nu_{\btheta} = \mbox{dim}(\btheta)$ et $p(\nu_{\btheta},n-\nu_{\btheta}; \alpha) = p(F(\nu_{\btheta},n-\nu_{\btheta}) \geq \alpha) $, $F$ étant le $100$ $(1-\alpha)$ percentile d'une variable de Fisher à $\nu_{\btheta}$ et $n-\nu_{\btheta}$  degrés de liberté \citep{Tomassone_etal}.

La fonction
\begin{equation}
f(t;\hat{\btheta}) =  \sum_{k=1}^K \pi_k(t;\hat{\bsw}) \hat{\bsbeta}_k^T\bt
\end{equation}
représente la courbe ajustée et 
\begin{equation}
s^2(t;\hat{\bsPhi}) = \frac{1}{n}{D(\hat{\bsPhi})}^T {I(\hat{\bsPhi})}^{-1} D (\hat{\bsPhi}) 
\label{eq. variance courbe ajustee}
\end{equation}
sa variance empirique, %
\begin{equation}
I(\hat{\bsPhi}) = -E \Big[\frac{\partial^2 \log p(y|t;\bsPhi)}{\partial \bsPhi \partial \bsPhi^T}\Big]_{\bsPhi=\hat{\bsPhi}}
\end{equation}
étant la matrice d'information de Fisher et
\begin{equation}
D(\hat{\bsPhi}) = \frac{\partial  f(t;\btheta)}{\partial \bsPhi}\!\Big \vert_{\bsPhi=\hat{\bsPhi}}
\end{equation}
le gradient de $f(t;\btheta)$  \citep{Jiang_and_tanner_IEEEinft_99}, calculés avec le paramètre estimé $\hat{\bsPhi}$.

Dans le cas non-linéaire, quand $n-\nu_{\btheta}$ est grand, (\ref{eq. region de confiance cas lineaire}) peut \^etre approximée asymptotiquement  \citep{Tomassone_etal, Gauchi_vila_coroller} et on obtient ainsi la région de confiance approchée au niveau $1-\alpha$ suivante :

\begin{equation}
f(t;\hat{\btheta}) \pm  \sqrt{\chi^2_{\nu_{\btheta},1-\alpha}} s (t;\hat{\bsPhi}).
 \end{equation}
\section{Expérimentation sur des données simulées}

L'objet de cette partie est d'évaluer l'approche de régression proposée en utilisant des données simulées. Pour ce faire, nous la
comparons à deux méthodes :
\begin{itemize}
\item une méthode de régression polynomiale par morceaux \citep{McGee}, \citep{sameSFC2007}, \citep{chamroukhi_et_al_ESANN2009} où l'estimation des paramètres est effectuée par une procédure de programmation dynamique \citep{bellman, stone},
\item une méthode de régression à processus markovien caché, c'est-à-dire où le processus caché $(z_1,\ldots,z_n)$ est supposé être une chaîne de Markov \citep{fridman}, \citep{rabiner}. Les paramètres de ce modèle sont estimés par l'algorithme de Baum-Welch \citep{BaumWelch}.
\end{itemize}
\medskip

 La qualité des estimations fournies par chacune des deux
méthodes est l'écart quadratique moyen (EQM) $\frac{1}{n}\sum_{i=1}^n (f_{sim}(t_i) - f_{est}(t_i))^2$
 entre la courbe de
régression estimée et la courbe de régression simulée où :
\begin{itemize}
\item $f_{est}(t_i)=\sum_{k=1}^{K} \pi_{k}(t_i;\hat{\bsw}) \hat{\bsbeta}^T_k \bt_{i}$  pour le modèle proposé,
\item $f_{est}(t_i)=\sum_{k=1}^{K}\hat{z}_{ik}\hat{\bsbeta}^T_{k}\bt_{i}$ pour le modèle de régression par morceaux,
\item $f_{est}(t_i)=\sum_{k=1}^{K}\omega_{k}(t_i;\hat{\bsPhi})\hat{\bsbeta}^T_k \bt_{i}$ pour le modèle markovien où les $\omega_{k}(t_i;\hat{\bsPhi})$ sont les probabilités dites de filtrage $\omega_{k}(t_i;\hat{\bsPhi}) = p(z_i=k|x_1,\ldots,x_i;\hat{\bsPhi})$ qui se calculent via une procédure du type ``forward-backward'' \citep{rabiner}.
\end{itemize}
\medskip
Chaque jeu de données est généré en ajoutant un bruit gaussien
centré à des points d'une courbe non linéaire, régulièrement
échantillonnés sur l'intervalle temporel $[0;5]$.

\subsection{Paramètres de simulation et réglage des algorithmes}

\begin{figure}[!h]
\centering \begin{tabular}{c}situation 1\\
\includegraphics[width=8.5cm,height = 4cm]{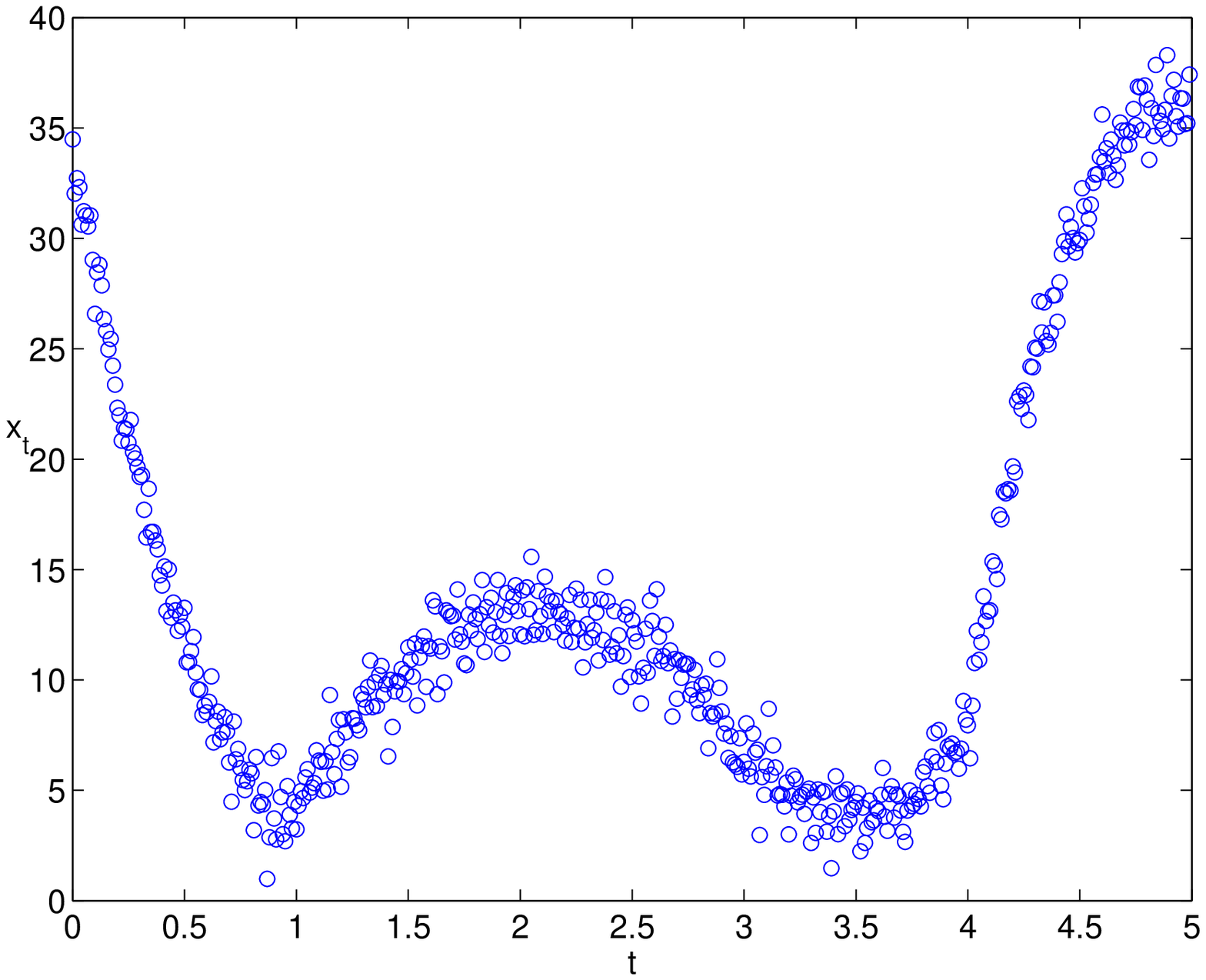}\\situation 2\\
\includegraphics[width=8.5cm,height = 4cm]{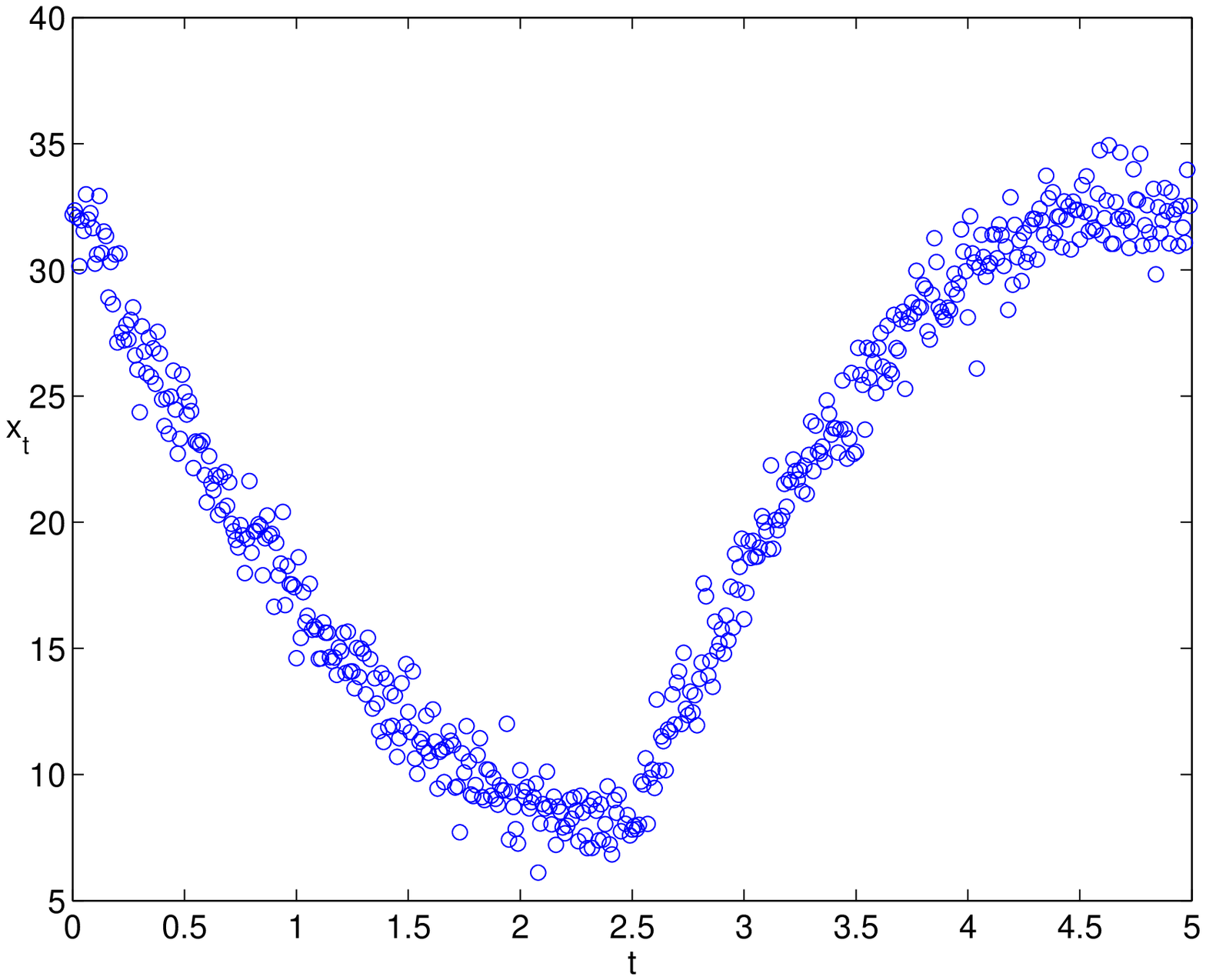}\\situation 3\\
\hspace*{-.2cm}\includegraphics[width=8.6cm,height = 4cm]{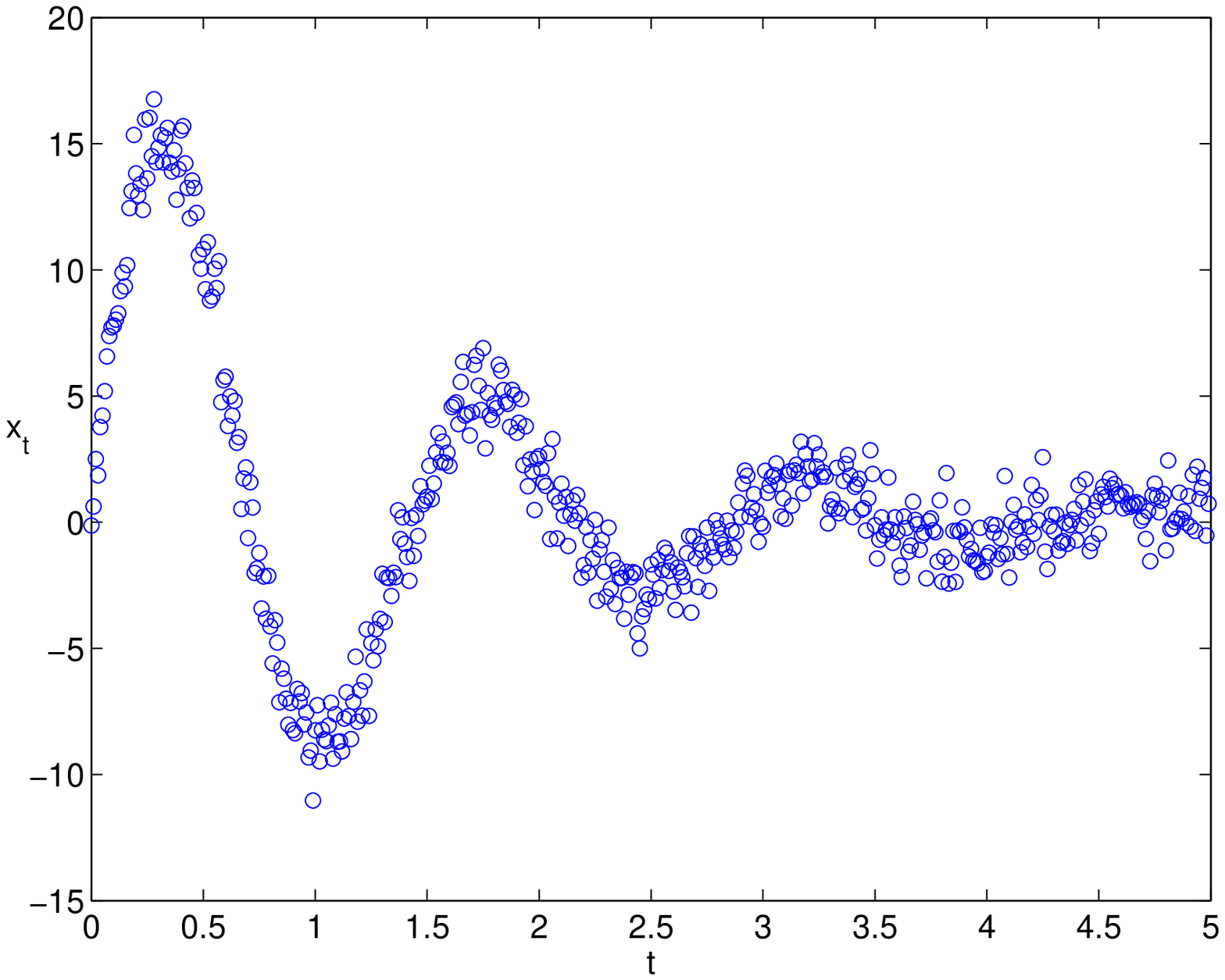}
  \end{tabular}
    \caption{Exemple de données simulées à partir des trois courbes}
      \label{signaux_sim}
\end{figure}

Trois fonctions non linéaires de régression ont été considérées. Ces
fonctions et leurs paramètres associés sont fournis dans le tableau
\ref{param}. Des exemples de données simulées à partir de ces
courbes sont représentés sur la figure \ref{signaux_sim}.
\begin{table}[!h]
\centering
{\small
\begin{tabular}{|l|ll|}
\hline
Fonctions & Paramètres&\\\hline
$f_1(t,\btheta)=\sum_{k=1}^4\pi_k(t;\bsw)\bsbeta_k^T \bt$&$\bsbeta_1=[34,-60,30]$ & $\bsw_1=[547,-154]$ \\
&$\bsbeta_2=[-17,29,-7]$ & $\bsw_2=[526,-135]$ \\
&$\bsbeta_3=[185,-104,15]$ & $\bsw_3=[464,-115]$ \\
&$\bsbeta_4=[-804,343,-35]$ & $\bsw_4=[0,0]$ \\ \hline
$f_2(t,\btheta)=\bsbeta_1^T \bt \mathds{1}_{[0 ; 2.5]}(t) + \bsbeta_2^T \bt \mathds{1}_{]2.5 ;
5]}(t)$&$\bsbeta_1=[33,-20,4]$  &\\
    &$\bsbeta_2=[-78,47,-5]$  &\\
\hline
$f_3(t)=20\sin(1.6\pi t)\exp(-0.7 t)$ & & \\
\hline
\end{tabular}}
\caption{Expressions analytiques des courbes de régression utilisées avec leurs paramètres}
\label{param}
\end{table}

Les trois algorithmes testés ont été lancés avec $(K=4,p=2)$ pour la
situation 1, \linebreak $(K=2,p=2)$ pour la situation 2 et avec
$(K=5,p=3)$ pour la situation 3. Pour chaque jeu de données,
l'algorithme EM a été lancé à partir de 10 initialisations
aléatoires différentes et seule la solution ayant la plus grande
vraisemblance a été retenue.

\medskip
\medskip

Il faut souligner ici que la première situation a tendance à favoriser l'approche proposée mais permet de la valider. Les deux autres fonctions considérées sont plus pertinentes pour évaluer la méthode proposée. Notamment, la troisième situation constitue un exemple typique  de fonction non linéaire.

\subsection{Résultats}

\begin{figure}[!h]
\centering
\begin{tabular}{cc}
\includegraphics[width=6.2cm,height=4.9cm]{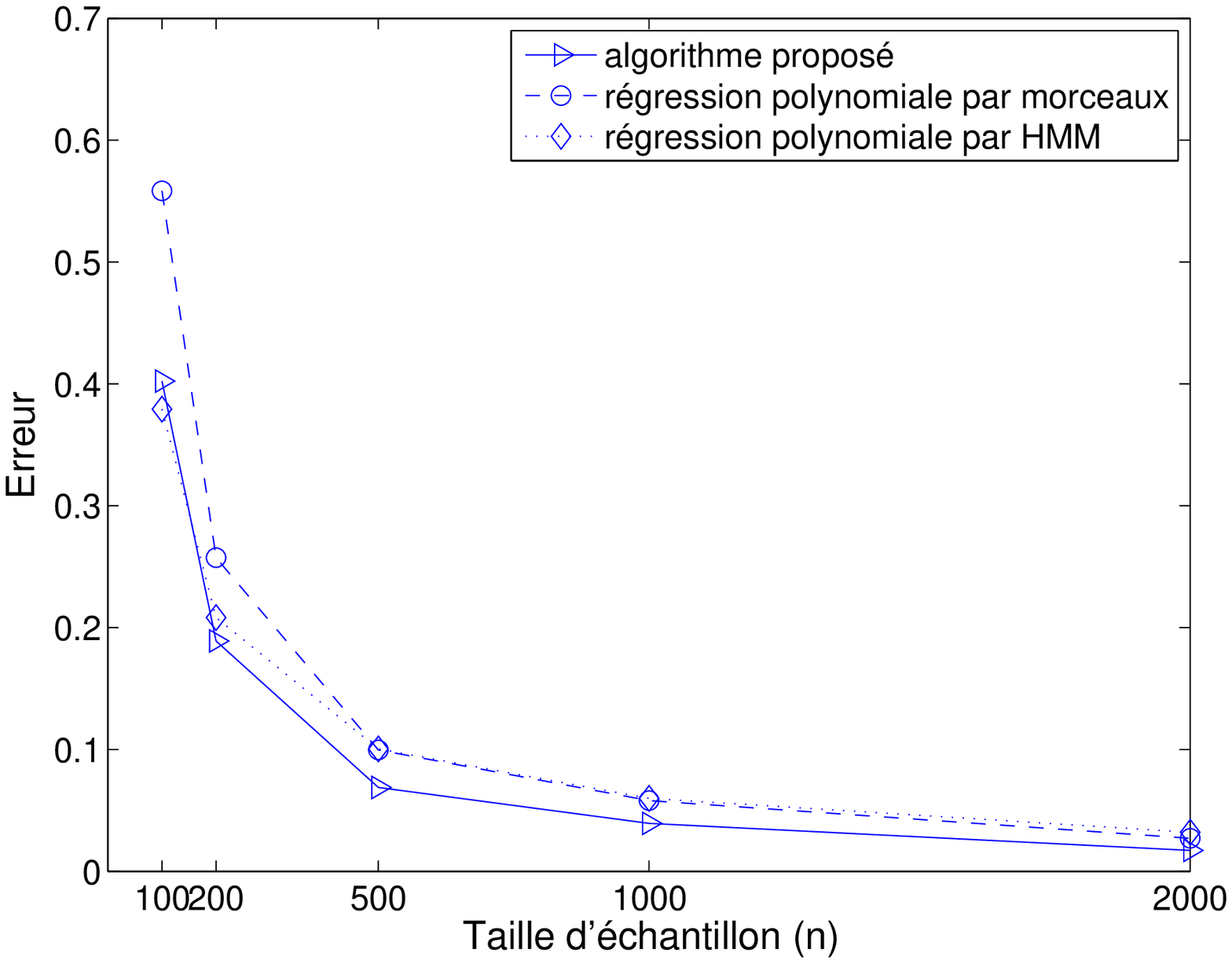}&
\includegraphics[width=6.2cm,height=4.9cm]{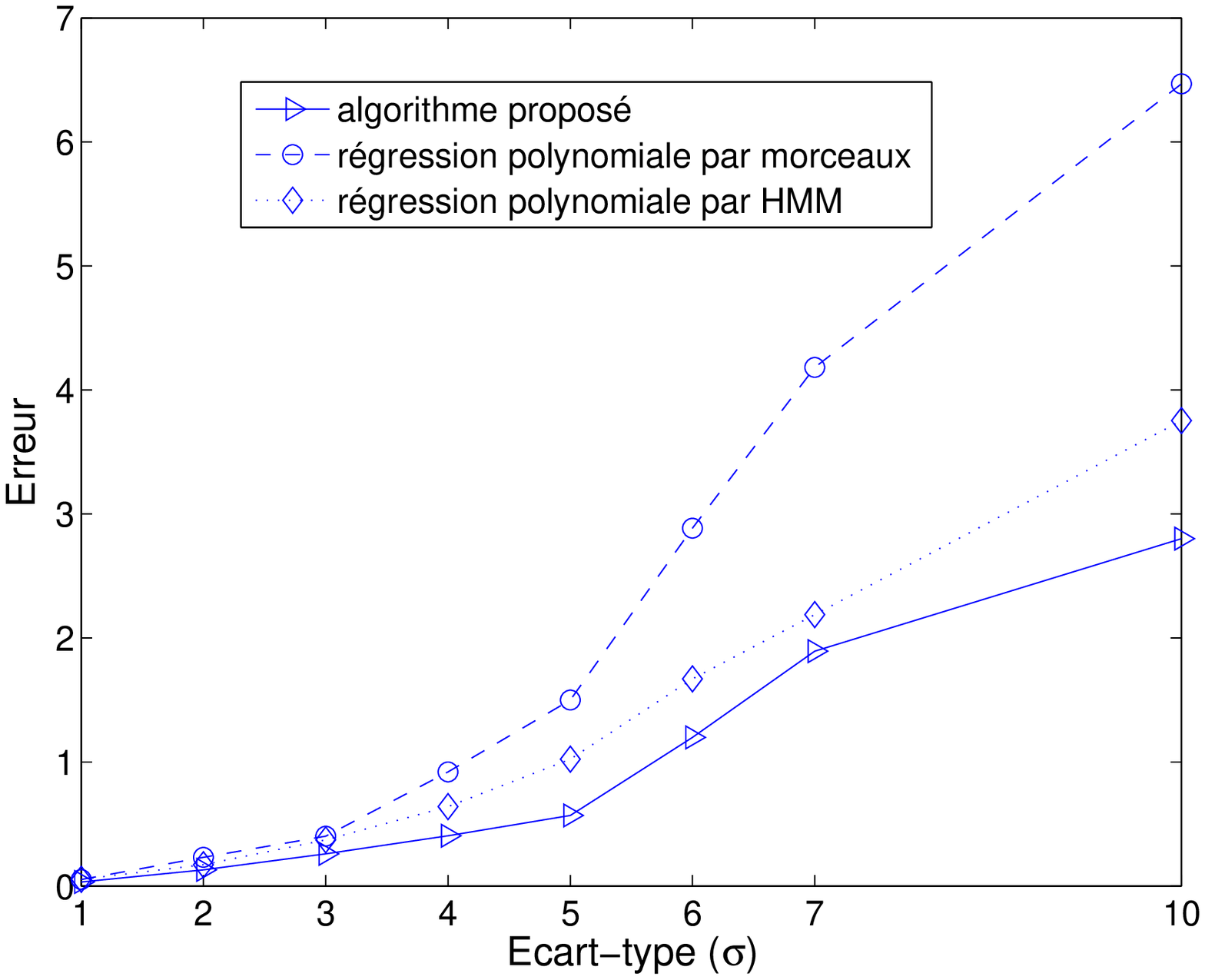}\\
\includegraphics[width=6.2cm,height=4.9cm]{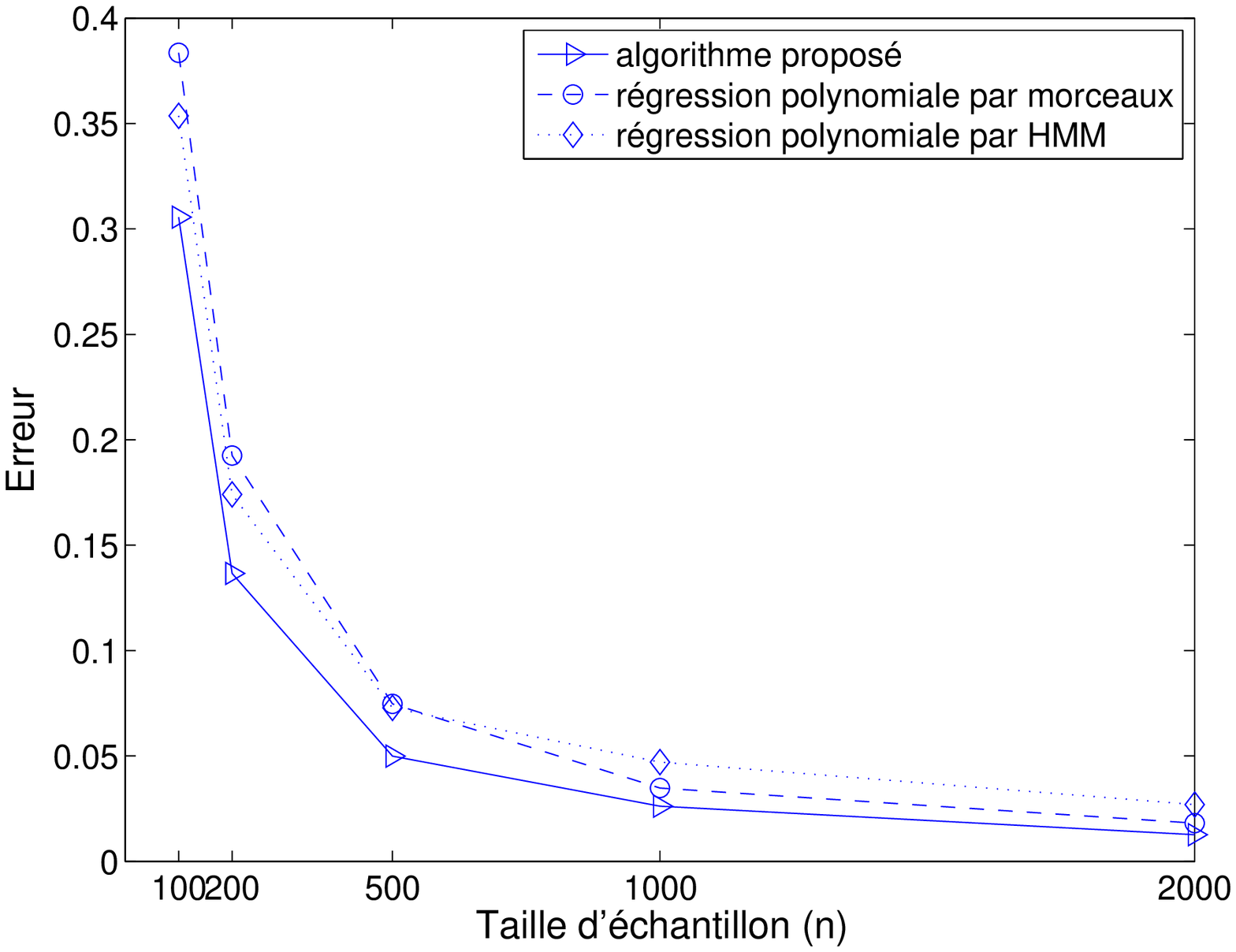}&
\includegraphics[width=6.2cm,height=4.9cm]{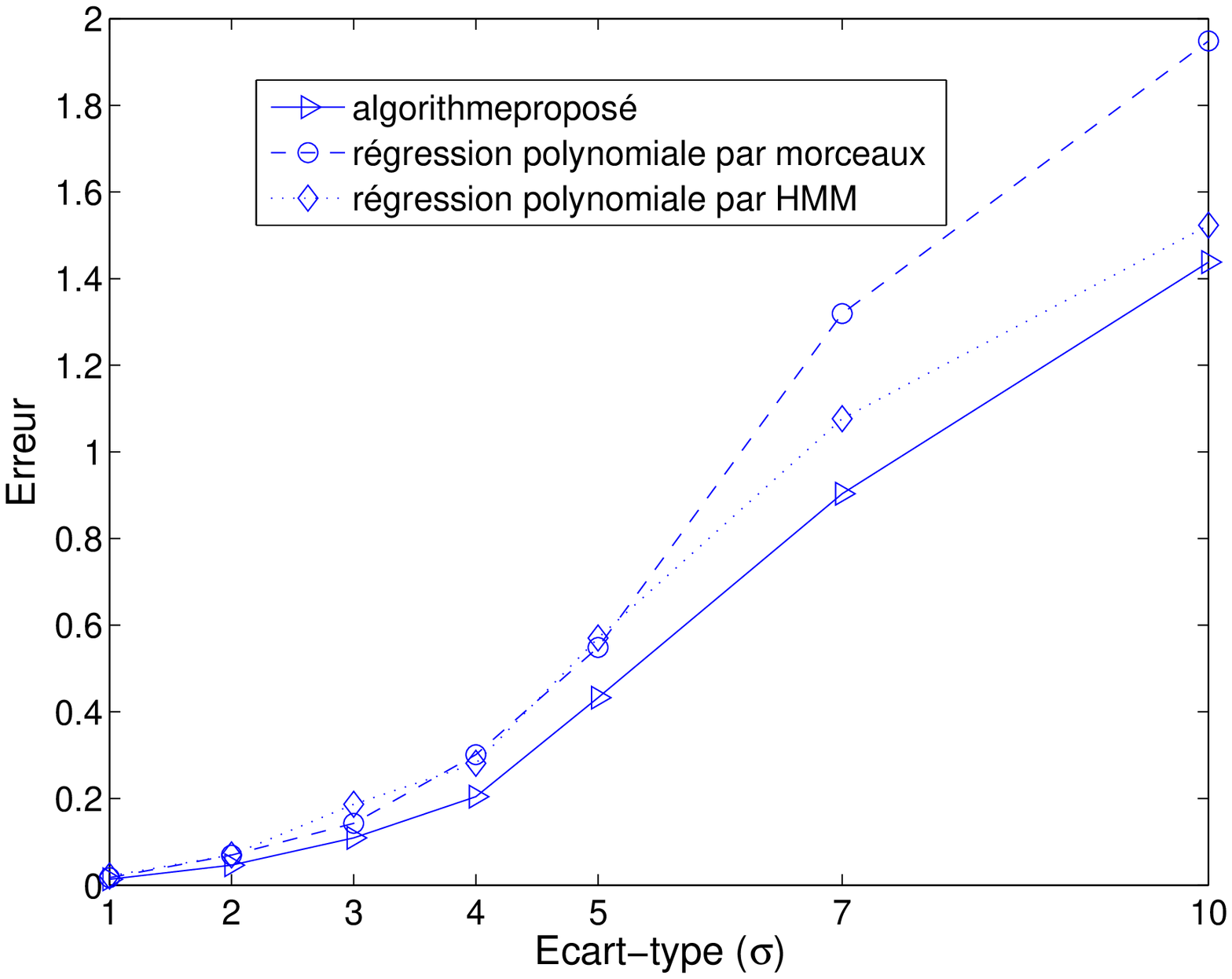}\\
\includegraphics[width=6.2cm,height=4.9cm]{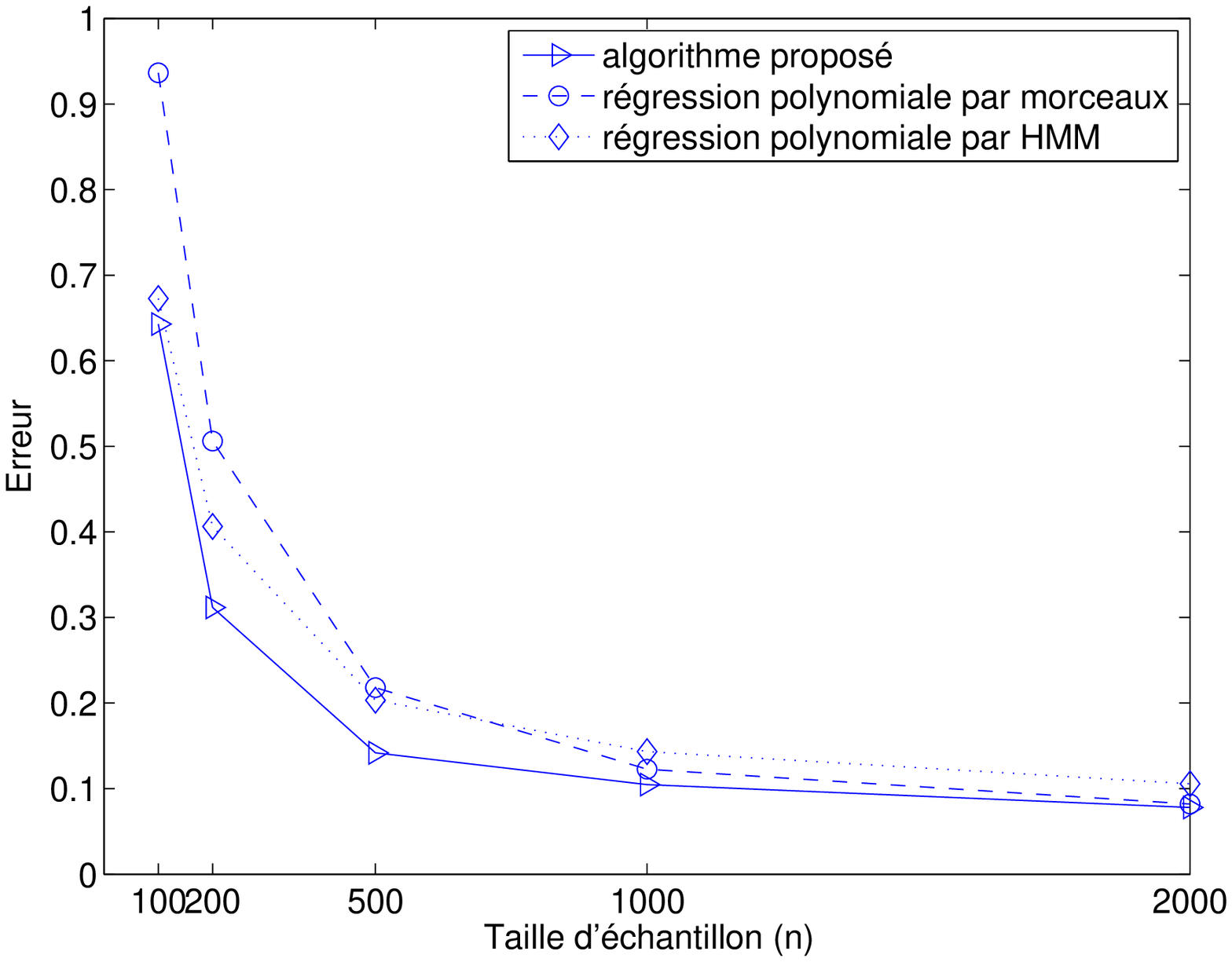}&
\includegraphics[width=6.2cm,height=4.9cm]{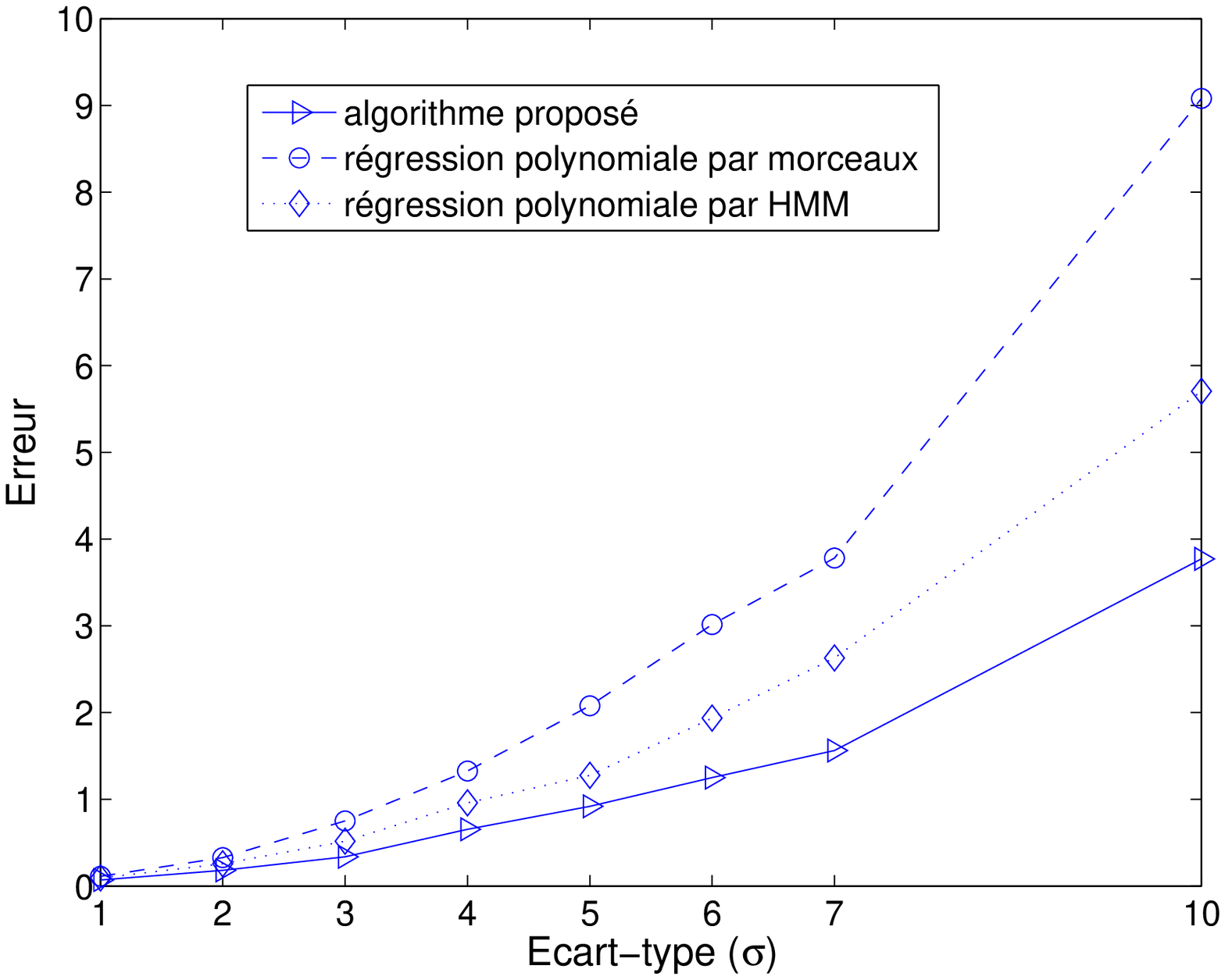}\\
\end{tabular}
\caption{Erreur entre la courbe estimée et la courbe réelle simulée
en fonction de la taille d'échantillon (gauche) ($\sigma=1.5$) et en
fonction de la variance du bruit (droite) ($n=500$), pour la
situation 1 (haut), la situation 2 (milieu) et  la situation 3
(bas)}\label{situation 1}
\end{figure}

La figure \ref{situation 1} montre pour chacune des situations,
comment varie l'écart entre les courbes de régression simulées et
les courbes estimées, en fonction de la taille d'échantillon et de
la variance du bruit. Pour chaque taille d'échantillon et chaque
valeur de la variance du bruit, l'erreur quadratique présentée
correspond à une moyenne sur 20 jeux de données différents.
 
Ces graphiques montrent que la méthode proposée donne de meilleurs
résultats que la méthode de régression polynomiale par morceaux et
la méthode basée sur un processus markovien caché. Cette différence entre les résultats fournis par les trois méthodes s'explique par leur gestion différente des passages d'un sous-modèle à l'autre. En effet, contrairement aux deux alternatives qui sont adaptées à des  signaux avec changement brusques, l'approche proposée permet de s'adapter à la fois aux transitions souples et brusques gr\^ace à la flexibilité de la fonction logistique qui modélise le processus latent.

\begin{figure}[!h]
\centering
\begin{tabular}{cc}
\includegraphics[width=6.25cm, height=3.8cm]{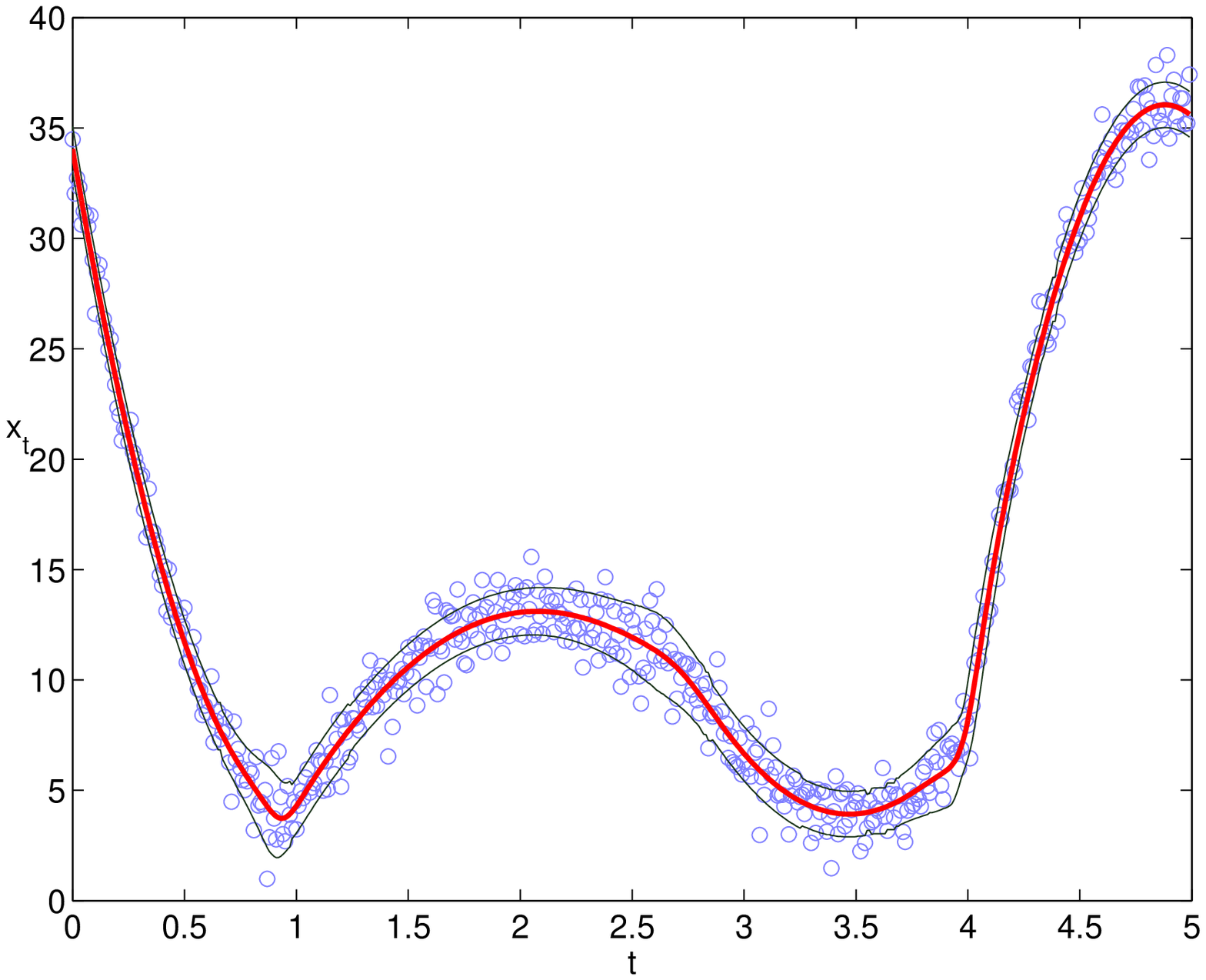}&\includegraphics[width=6.25cm,height=3.8cm]{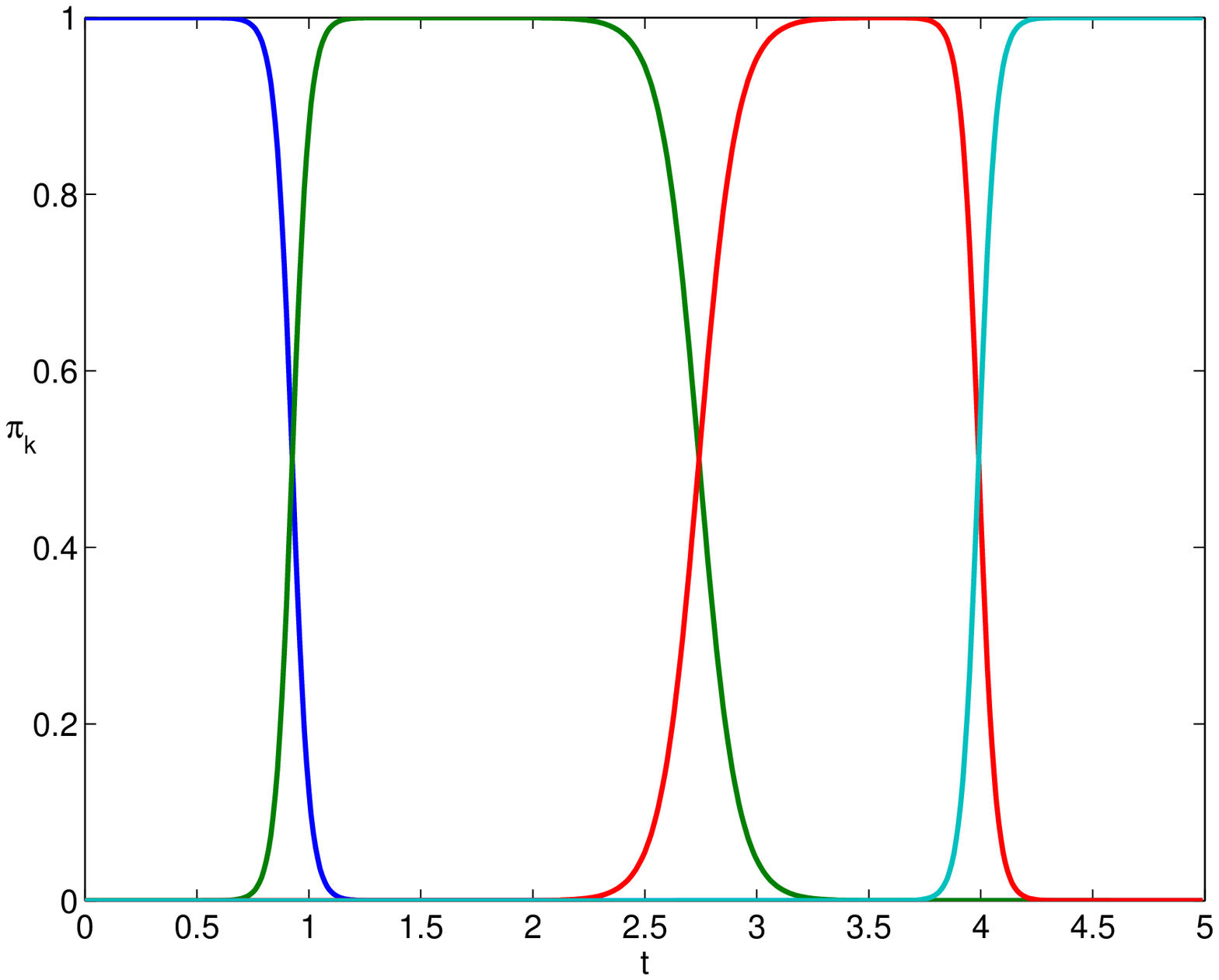}\\
\includegraphics[width=6.25cm, height=3.8cm]{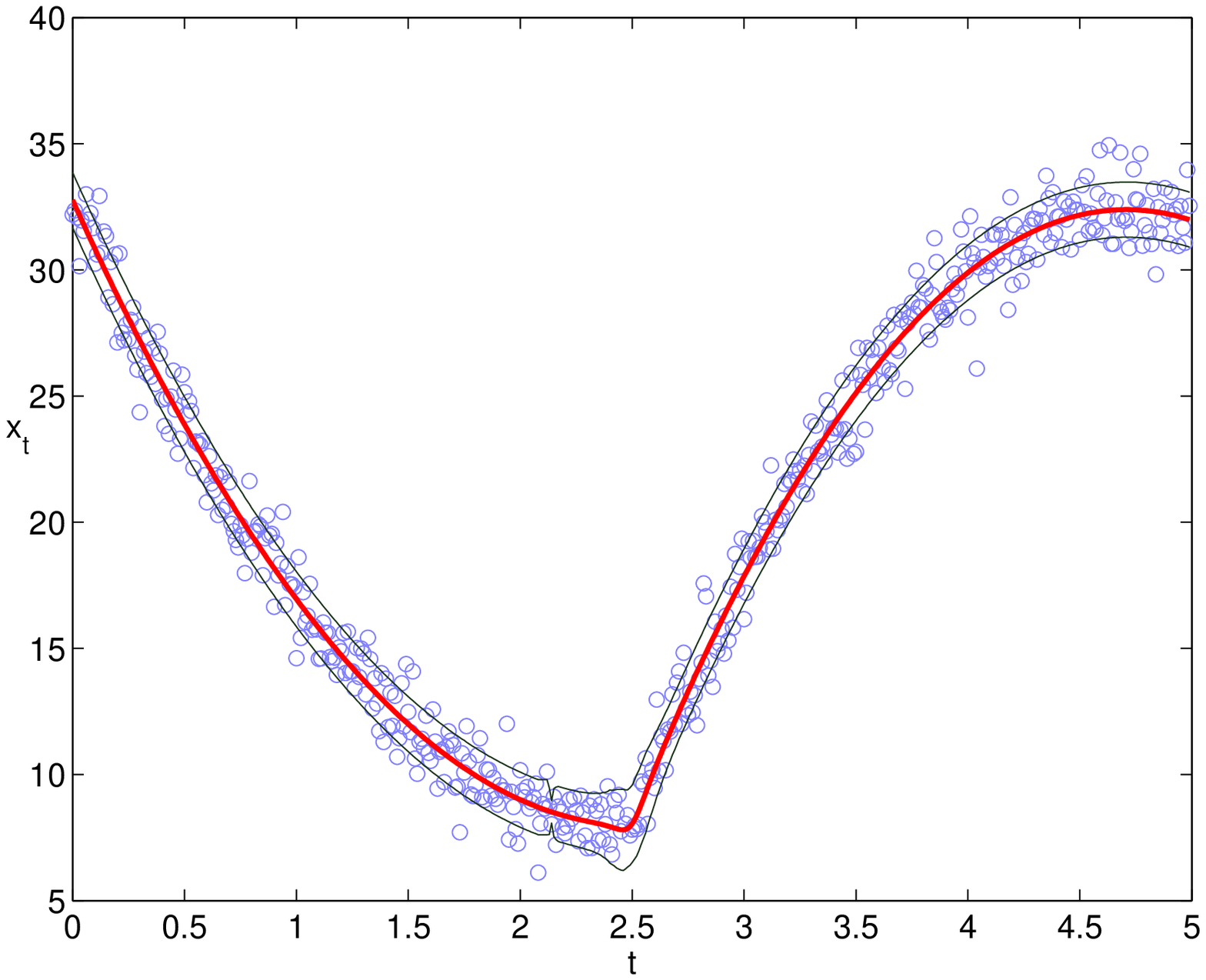}&\includegraphics[width=6.25cm,height=3.8cm]{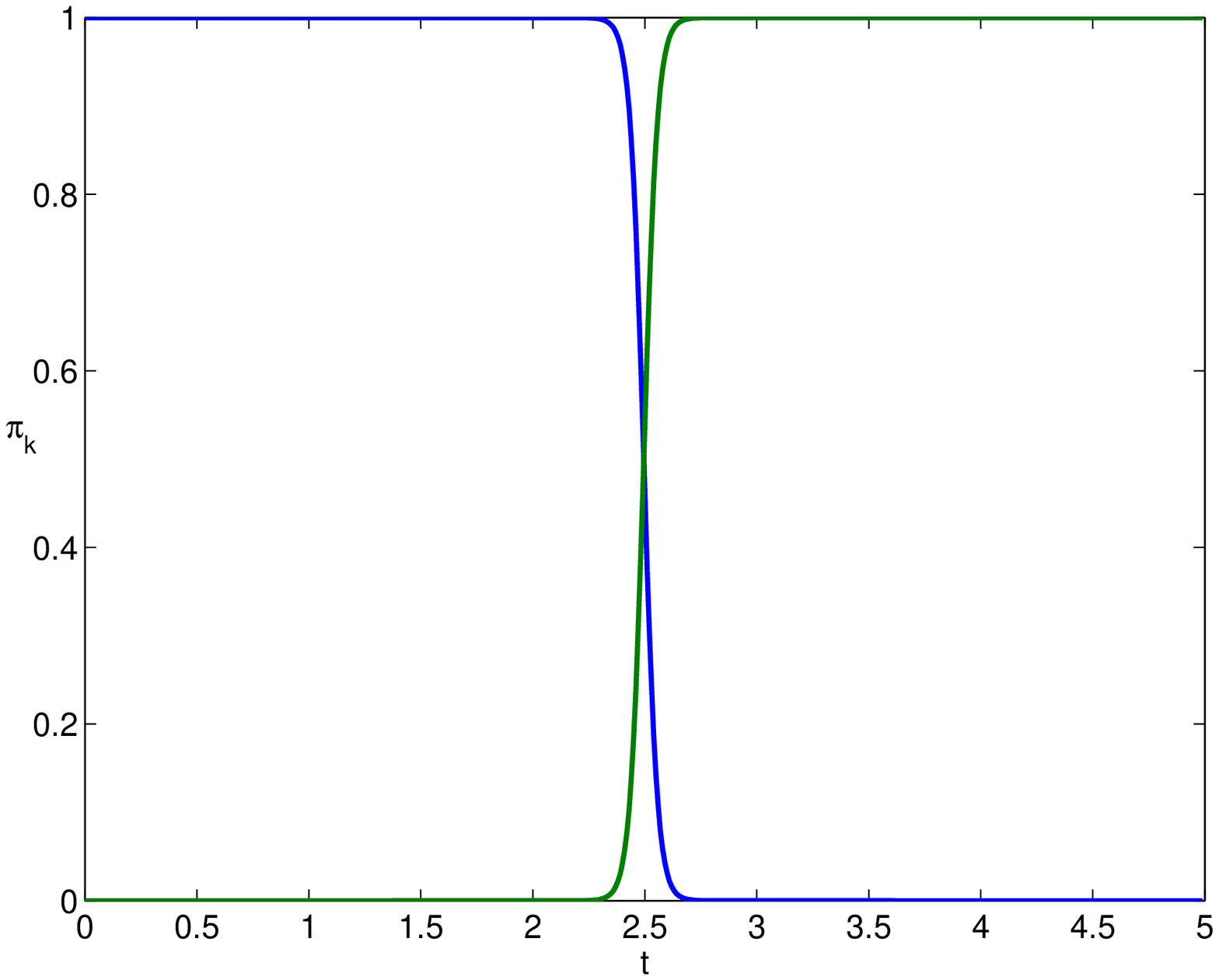}\\
\hspace*{-.2cm} \includegraphics[width=6.35cm, height=3.8cm]{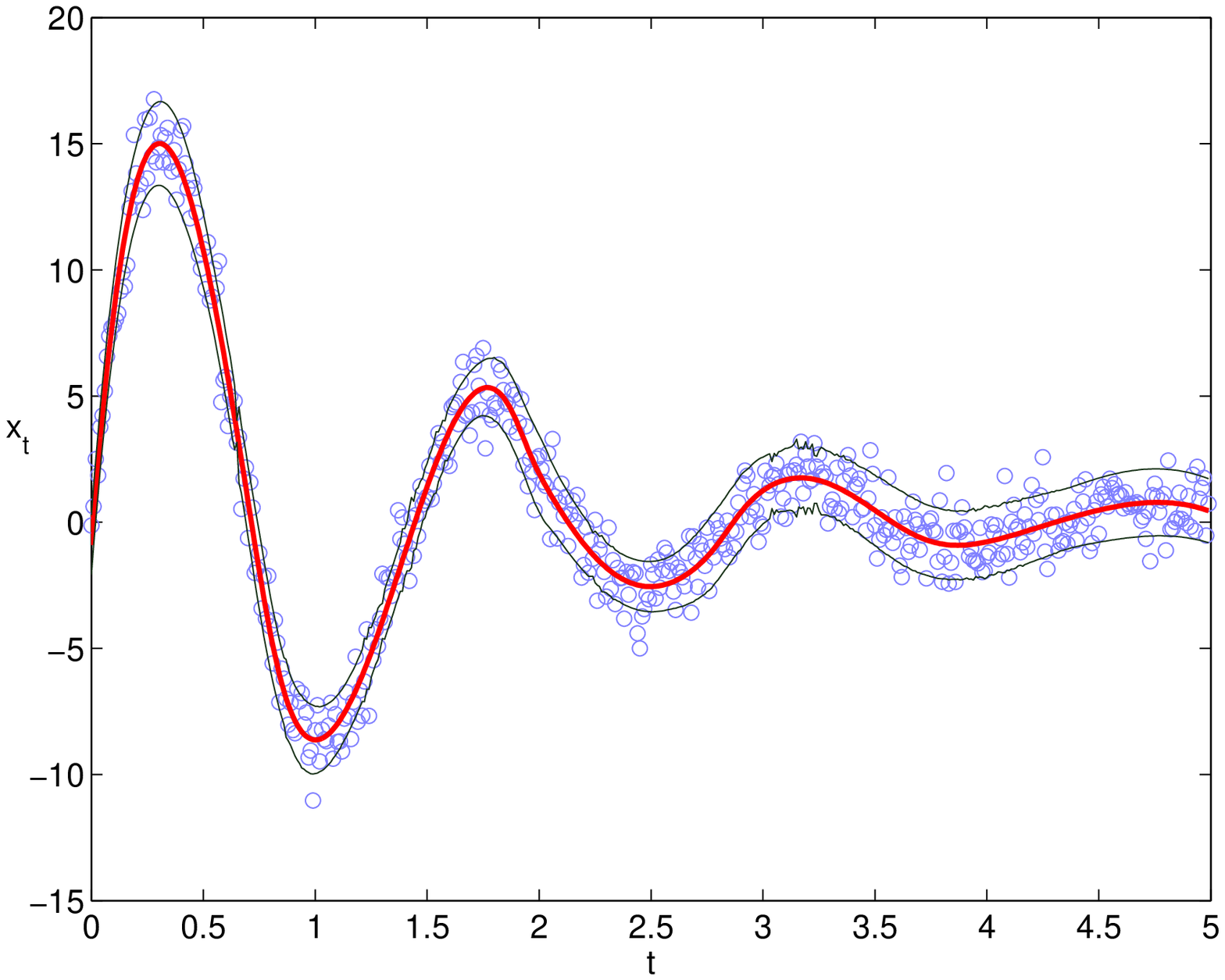}&\includegraphics[width=6.25cm,height=3.9cm]{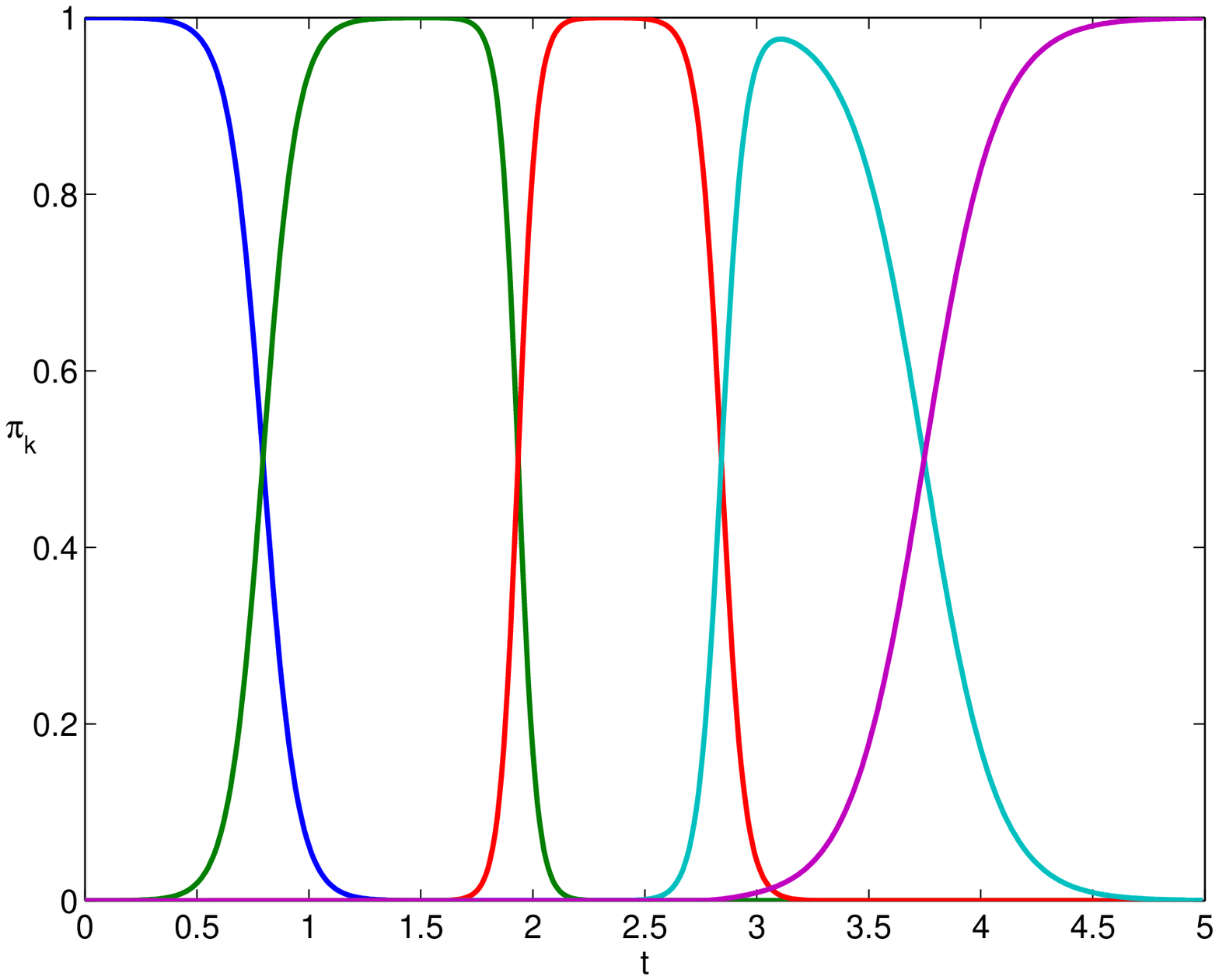}\\
\end{tabular}
\caption{Fonctions de régression $f(t;\btheta)$ estimées par
l'algorithme proposé et régions de confiance à 95 \% associées (à gauche) et proportions $\pi_k$
correspondantes (à droite)}
\label{estimations}
\end{figure}

En outre, on peut observer que l'écart entre les courbes estimées et les courbes
simulées décroît quand la taille d'échantillon augmente.
L'augmentation de la variance du bruit entraîne quant à elle une
augmentation de l'erreur qui est plus prononcée pour le modèle de
régression polynomial par morceaux. La figure \ref{estimations}
montre, pour chacune des situations, un exemple de courbe de
régression estimée avec l'algorithme proposé et la région de confiance à 95 \% correspondante.

Le tableau \ref{table. resultats BIC en fonction de K et p} montre les pourcentages moyens de choix  du nombre de sous-modèles de régression $K$  et du degré $p$ des polynômes, sélectionnés par le critère BIC, sur les signaux simulés de la situation 1. On observe que les vraies valeurs ($K=4$ et $p=2$) ont été sélectionnées avec un pourcentage maximum de 63 \%.
\begin{table}[!h]
\centering
{\small \begin{tabular}{l c c c c c c c}
\hline & $p$  &1&2&3&4&5&6\\
$K$ & & & & & & & \\
\hline
$2$& &     $0$    &  $0$                &    $0$  &   $13$    &  $3$  &    $0$   \\
$3$& &     $0$    &  $0$                &    $7$  &  $13$    &  $0$  &    $0$   \\
$4$& &     $0$    &  $\boldsymbol{63}$  &    $1$  &   $0$    &  $0$  &    $0$   \\
$5$& &     $0$    &  $0$                &    $0$  &   $0$    &  $0$  &    $0$    \\
$6$& &     $0$    &  $0$                &    $0$  &   $0$    &  $0$  &    $0$\\
$7$& &     $0$    &  $0$                &    $0$  &   $0$    &  $0$  &    $0$\\
\hline
\end{tabular}}
\caption{ \small Pourcentage  de choix de $(K,p)$ par le critère BIC obtenu  avec les signaux de la situation 1}
\label{table. resultats BIC en fonction de K et p}
\end{table}

Le tableau \ref{table. temps de calcul} montre les temps de calculs moyens obtenus pour les trois méthodes. On remarque que le modèle proposé et le modèle markovien, qui sont basés sur une méthode itérative de type EM pour l'estimation des paramètres, ont des temps de calculs quasi identiques  qui sont peu sensibles à la taille d'échantillon. Cependant, la méthode de régression par morceaux est très co\^uteuse en temps de calcul à cause de la procédure de programmation dynamique dont le temps de calcul augmente considérablement avec la taille d'échantillon.
\begin{table}[!h]
{\scriptsize \begin{tabular}{llccccccccc}
\hline
 & & & Situation 1 & &  & Situation 2 &   &  & Situation 3 & \\
\hline
 & \!\!\!\!\!\! Algorithme : &  algo 1& algo 2 & algo 3 & algo 1& algo 2 & algo 3 & algo 1& algo 2 & algo 3 \\
$n$ & & & & & & & & & & \\
\hline
$100$  & & $0.27$  & $0.15$  &  $0.12$  &  $0.08$  & $0.14$  & $0.07$ &  $0.35$  &  $0.15$ & $0.17$ \\
$200$  & & $0.45$  & $0.67$  &  $0.19$  &  $0.09$  & $0.66$  & $0.10$ &  $0.48$  &  $0.70$ & $0.31$ \\
$500$  & & $1.08$  & $5.09$  &  $0.46$  &  $0.21$  & $5.05$  & $0.21$ &  $1.77$  & $5.33$  & $0.74$  \\
$1000$ & & $1.93$  & $25.86$ &  $0.64$  &  $0.27$  & $25.65$ & $0.34$ &  $3.06$  & $26.59$ & $1.01$  \\
$2000$ & & $2.80$  & $147.7$ & $1.39$   &  $0.47$  & $147.2$ & $0.65$ &  $5.03$  & $148.2$ & $2.17$  \\
\hline
\end{tabular}}
\caption{ {\small Temps de calculs moyens en secondes, en fonction de la taille d'échantillon $n$, obtenus avec les trois méthodes : approche proposée (algorithme 1), approche de régression par morceaux (algorithme 2) et approche de régression à processus markovien (algorithme 3),  pour les trois situations des données simulées}}
\label{table. temps de calcul}
\end{table}

\section{ Expérimentation sur des données réelles}

Dans le cadre d'une application de suivi d'état de fonctionnement du
mécanisme d'aiguillage des rails, nous avons été amenés à
paramétriser des signaux non linéaires représentant la puissance
consommée par le moteur d'aiguillage durant des man{\oe}uvres
d'aiguillage. Cette paramétrisation vise à représenter dans un
espace de dimension peu élevée ces signaux non linéaires qui sont formés de 562 points. Nous avons réalisé cette tache
par la méthode de régression à processus latent proposée, les
paramètres de régression étant directement utilisés comme paramètres
des signaux.

Une man{\oe}uvre d'aiguillage est constituée de mouvements mécaniques des différents organes
liés à l'aiguille, qui sont mobilisés successivement. Ces mouvements se traduisent, sur le
signal de puissance consommée en fonction du temps, par différentes phases de
fonctionnement (5 phases : Appel moteur, décalage-déverrouillage, translation, verrouillage-calage et friction). Le nombre des composantes régressives $K$ du modèle, pour chaque
signal considéré, a été donc fixé à $K=5$ qui correspond au nombre de
phases d'une man{\oe}uvre d'aiguillage. L'ordre du polynôme a été fixé à $p=3$
qui est adapté à la forme des signaux traités.

\begin{figure}[!h]
\centering \begin{tabular}{c}(a)\\
\includegraphics[width=7.5cm,height=4.4cm]{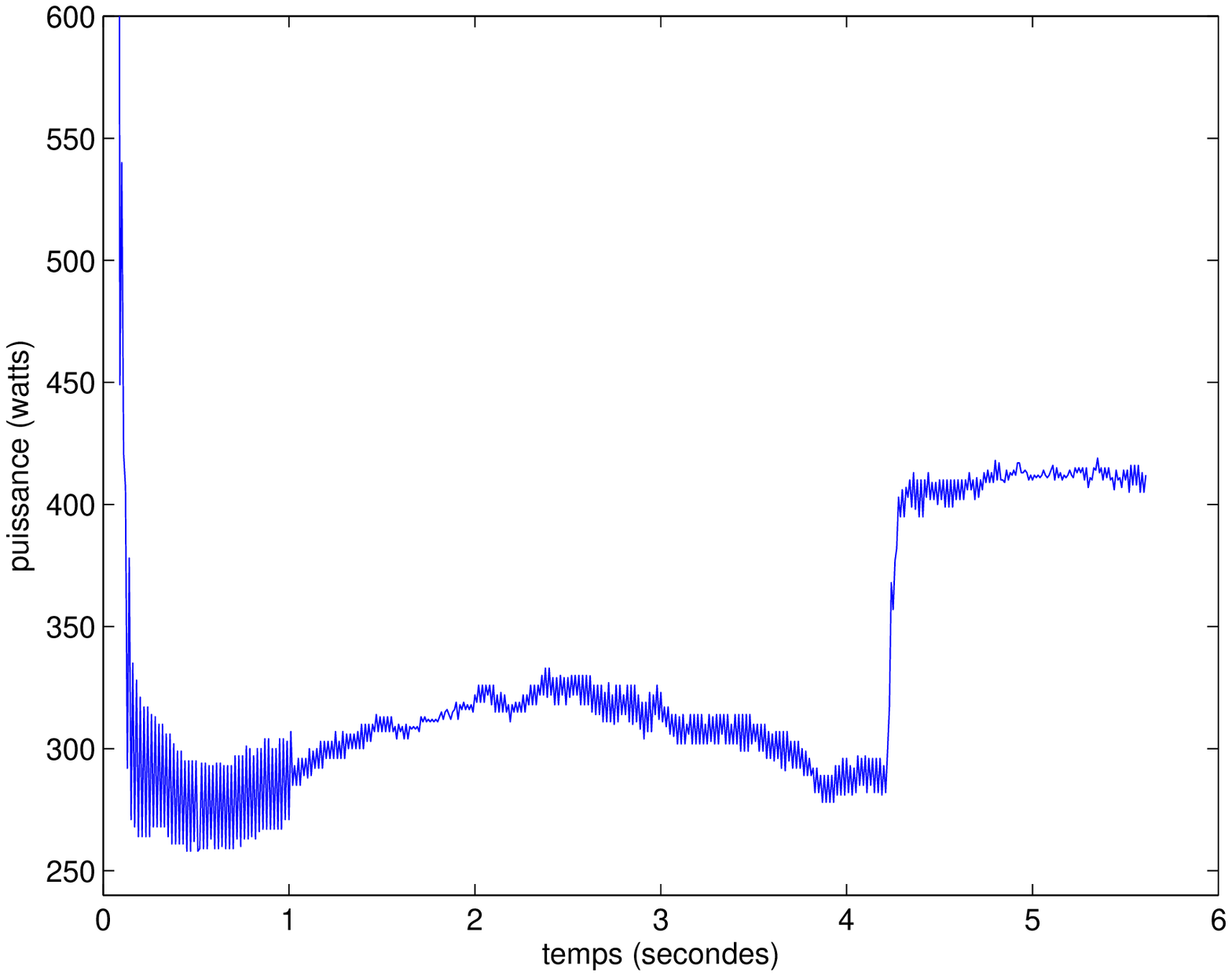}\\(b)\\
\includegraphics[width=7.5cm,height=4.4cm]{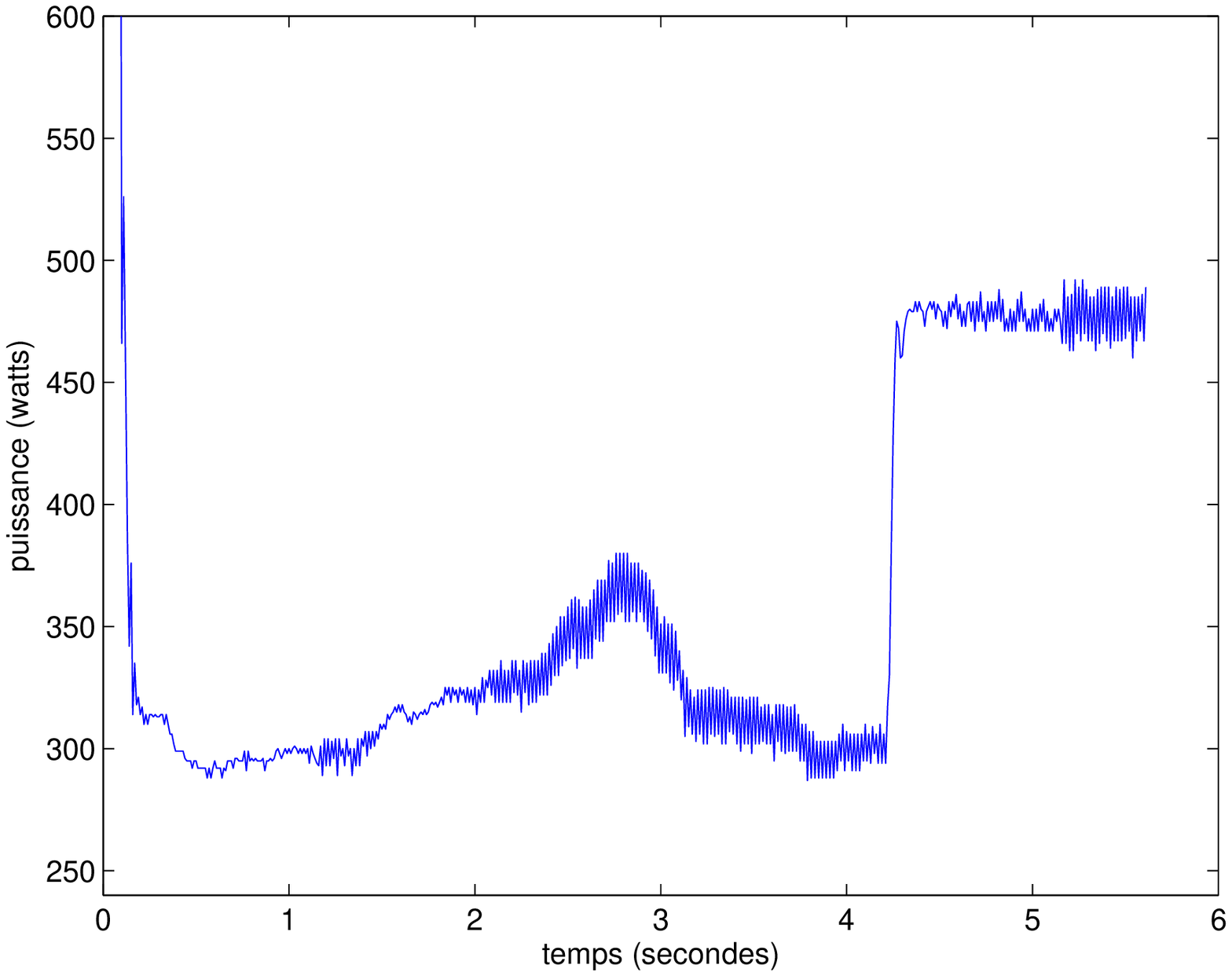}\\(c)\\
\includegraphics[width=7.5cm,height=4.4cm]{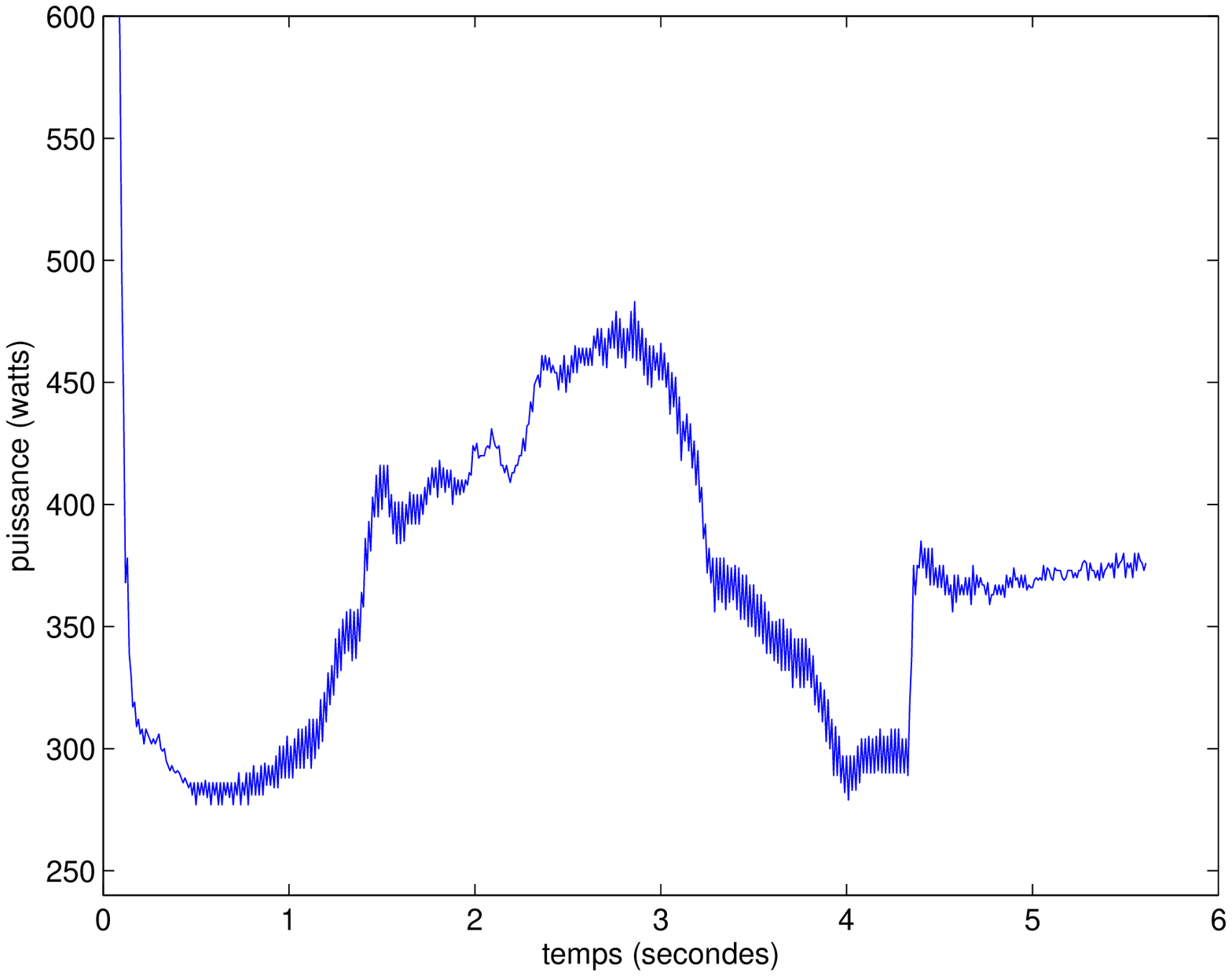}\\
  \end{tabular}
    \caption{Exemples de signaux de man{\oe}uvres d'aiguillage}
      \label{fig. exemple_signaux_reels}
\end{figure}

\medskip
La figure \ref{fig. exemple_signaux_reels} nous montre trois exemples
de signaux réels sur lesquels nous avons appliqué la méthode. Ces
signaux correspondent à trois états de fonctionnement : état sans
défaut (a), état avec défaut tolérable (b) et état avec défaut
critique (c).

\medskip

La qualité des estimations fournies par les différentes
méthodes est mesurée par l'écart quadratique moyen (EQM) entre le
vrai signal et son estimation. Ce critère est donné par $EQM =
\frac{1}{n}\sum_{i=1}^n (x_i - \hat{x}_i)^2$ .

Les résultats correspondants aux trois situations de signaux de
man{\oe}uvres d'aiguillages obtenus par les trois méthodes sont
donnés dans le tableau \ref{resultats MSE}. Ils mettent en évidence
les bonnes performances de la méthode proposée en termes d'écart
quadratique moyen. On remarque que, pour la situation (a), les écarts quadratiques moyens obtenus  par la régression par morceaux et la régression par HMM sont légèrement inférieurs à l'écart quadratique moyen obtenu avec la méthode proposé. Cela est du au fait que le signal de la situation (a) présente clairement des changements brusques entre les différentes phases.
\begin{table}[!h]
\centering
{\footnotesize \begin{tabular}{cccc}
\hline
Situation & Approche proposée & Régression par morceaux & Régression par HMM \\
\hline
(a) & $784.93$ & $\boldsymbol{781.64}$ & $783.20$\\
(b) & $\boldsymbol{1800.81}$ & $1928.92$ & $1816.31$\\
(c) & $\boldsymbol{309.80}$ & $310.25$ & $314.83$\\
\hline
\end{tabular}}
\caption{{\small  Erreurs quadratiques moyennes obtenues par les trois
méthodes pour les trois signaux de man{\oe}uvres d'aiguillage}}
\label{resultats MSE}
\end{table}

\begin{figure}[!h]
\centering
\begin{tabular}{cc}
\includegraphics[width=6.4cm,height=4.5cm]{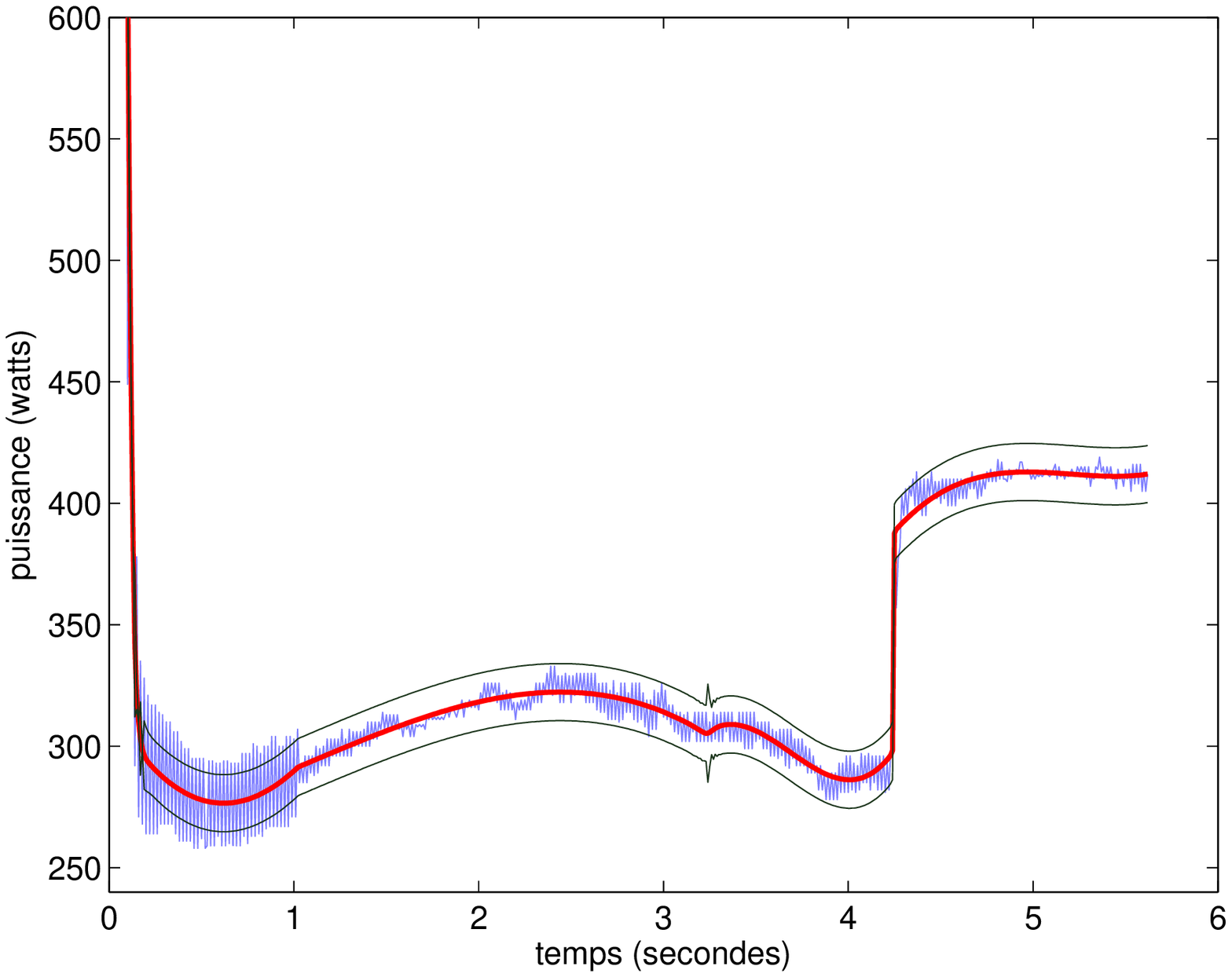}&\includegraphics[width=6.2cm,height=4.4cm]{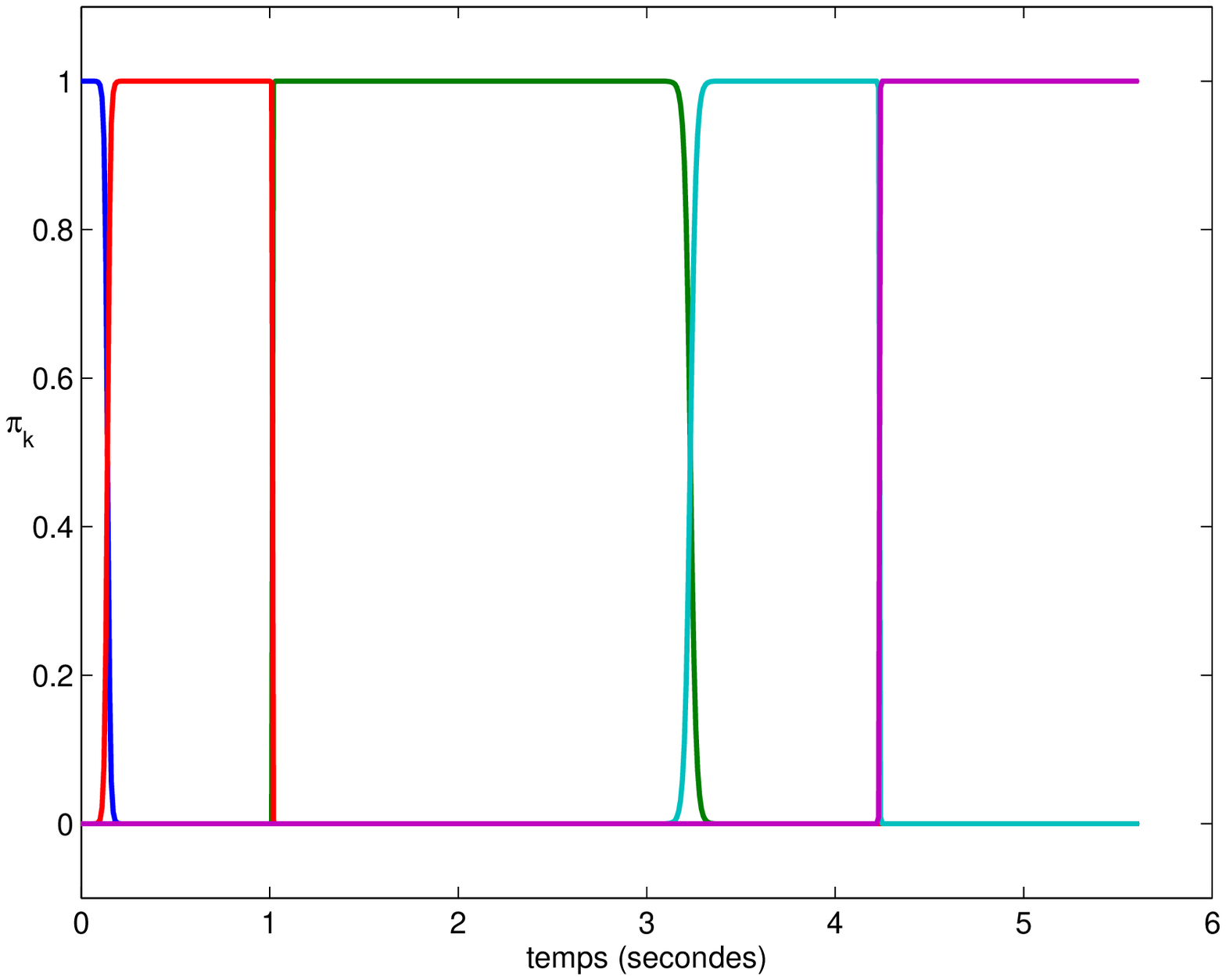}\\
\includegraphics[width=6.4cm,height=4.5cm]{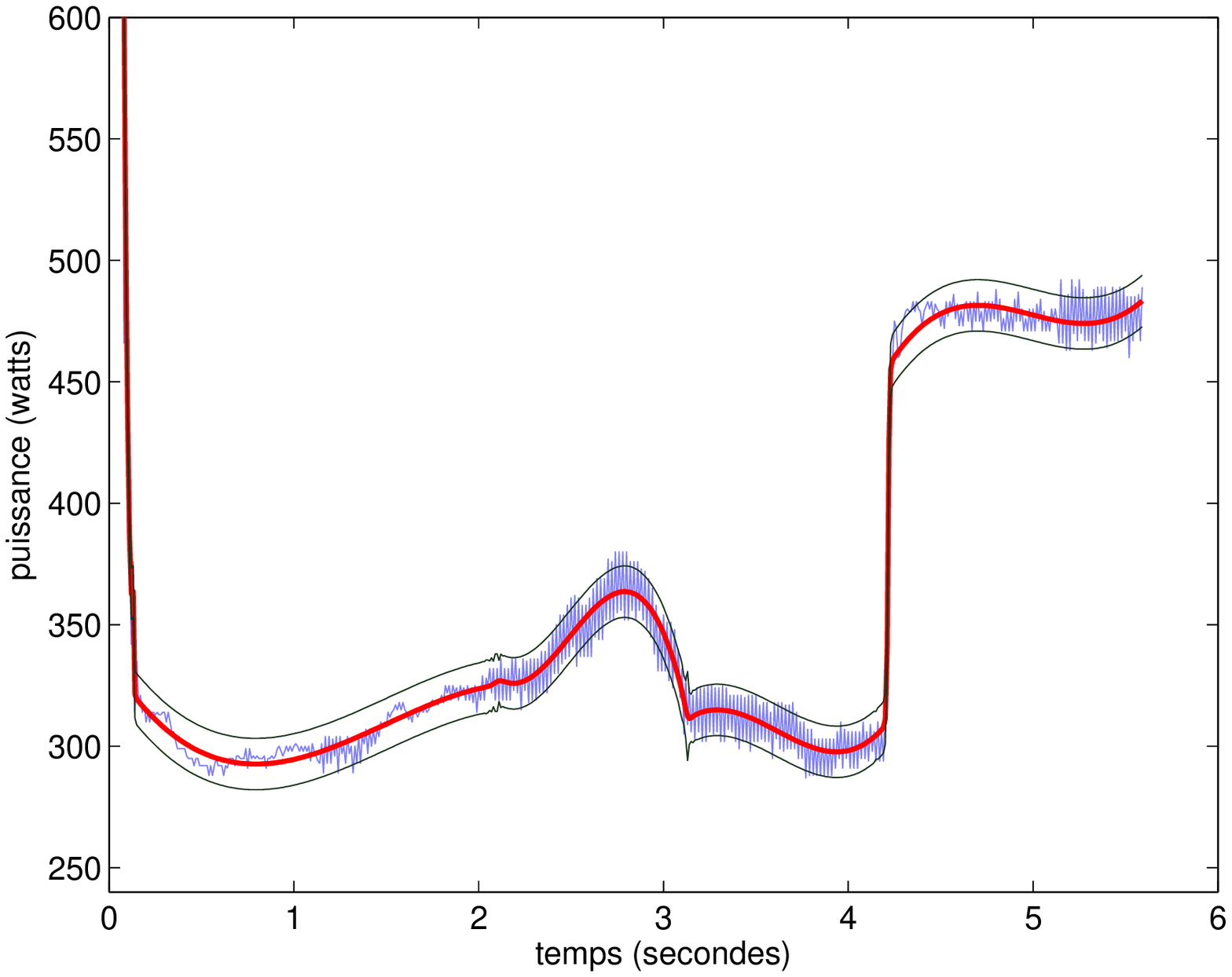}&\includegraphics[width=6.2cm,height=4.4cm]{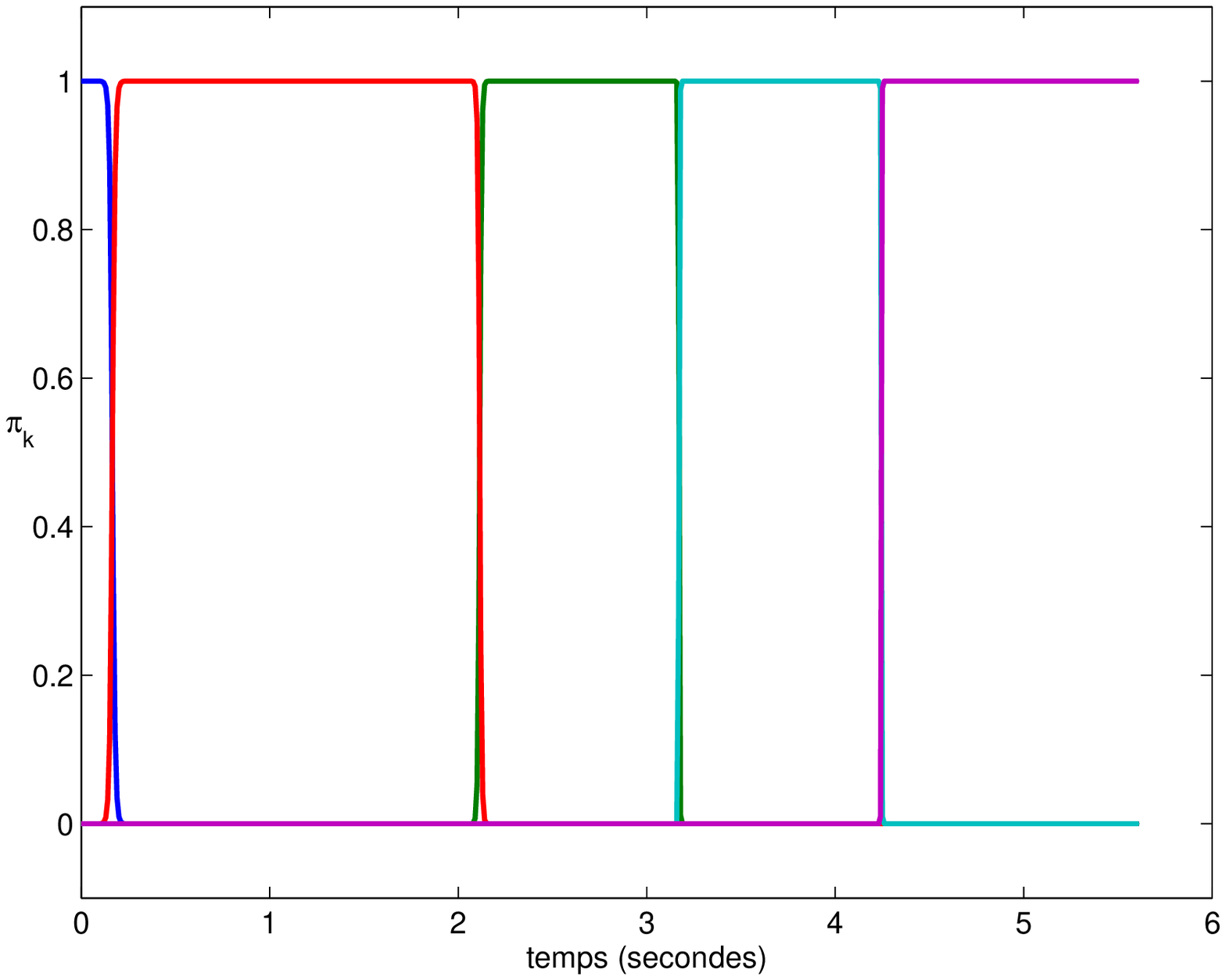}\\
\includegraphics[width=6.4cm,height=4.5cm]{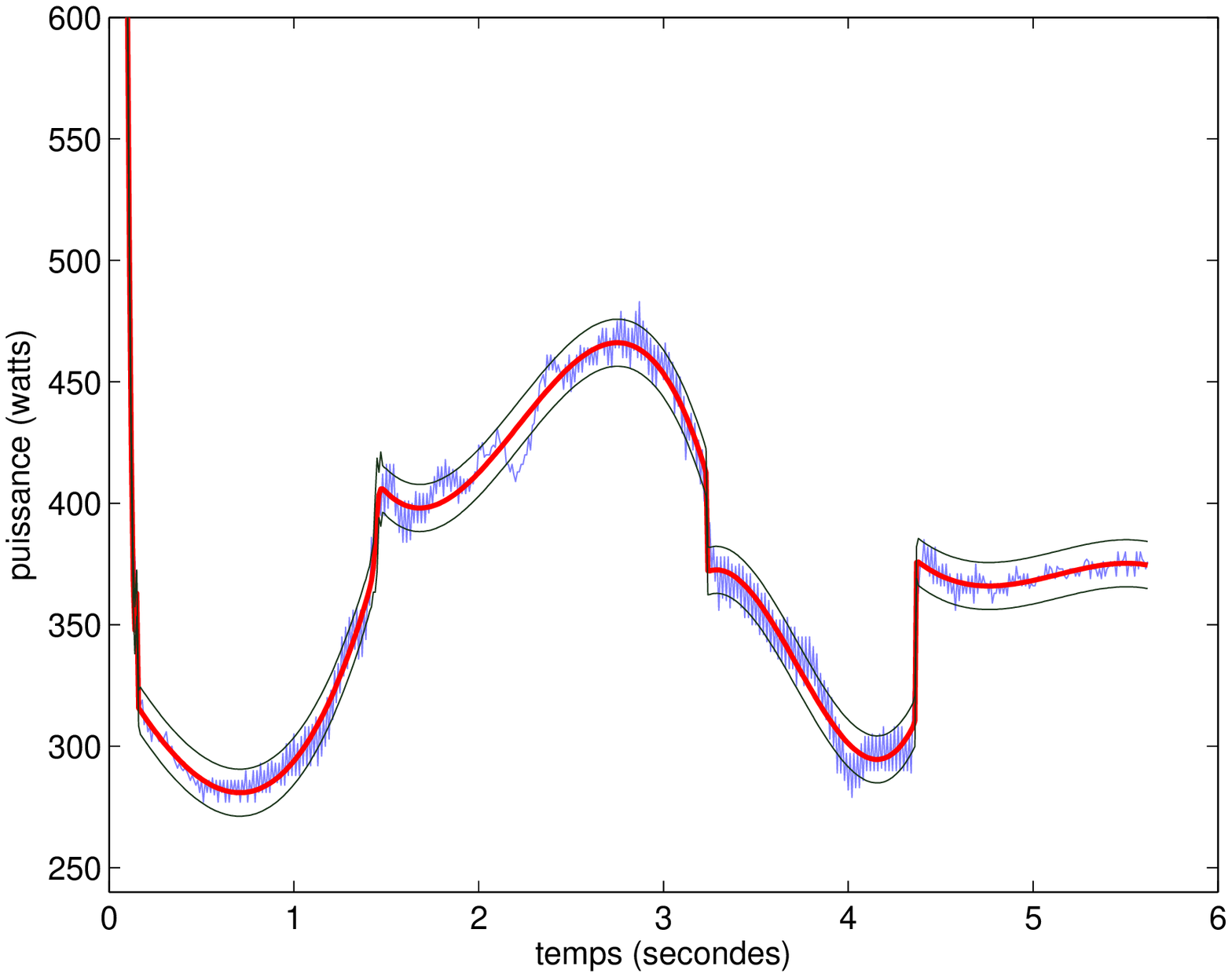}&\includegraphics[width=6.2cm,height=4.4cm]{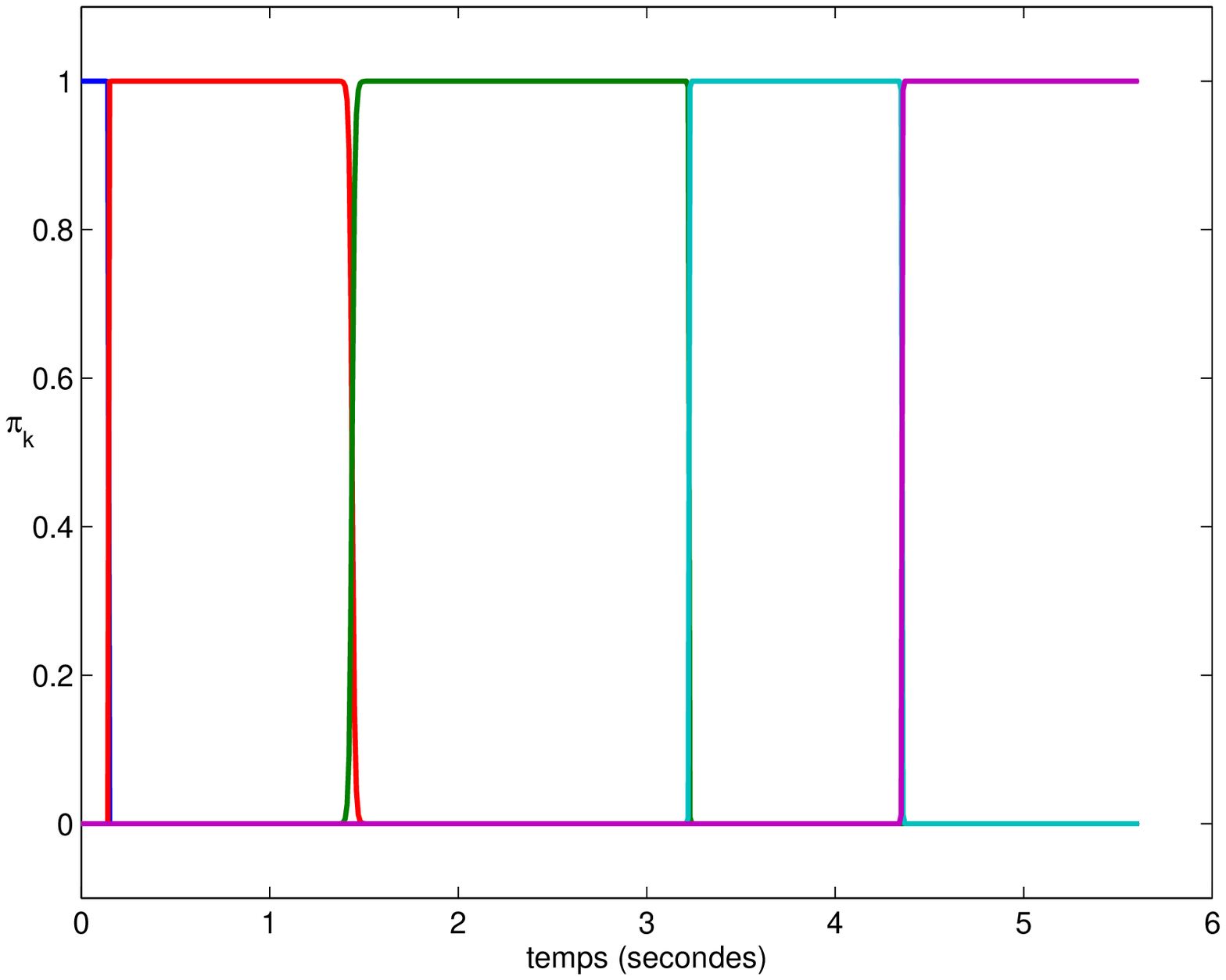}\\
\end{tabular}
\caption{Signal et modèle estimé par l'algorithme proposé et régions de confiance au niveau de confiance de 95 \% associées (gauche)
et proportions correspondantes (droite), pour la situation sans
défaut (haut), avec défaut mineur (milieu) et avec défaut critique
(bas) }\label{estimations_aig_sit1et2}
\end{figure}

La figure \ref{estimations_aig_sit1et2} montre les résultats
graphiques correspondant aux trois signaux présentés.  Les proportions du mélange estimées sont  cohérentes avec la réalité des phases des man{\oe}uvres d'aiguillage considérées.

\begin{figure}[!h]
\centering
\includegraphics[width=6.5cm]{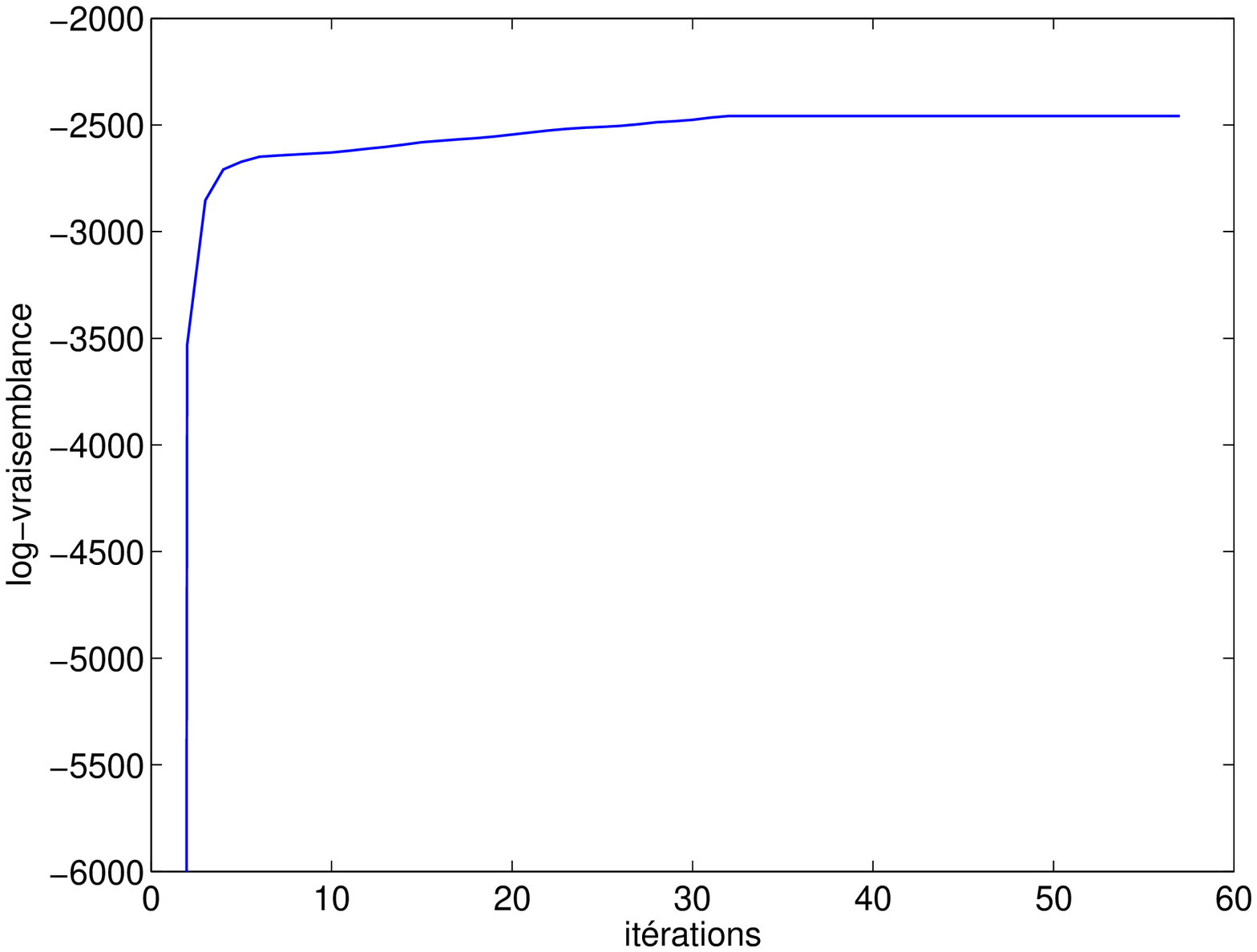}
\caption{Convergence du critère de log-vraisemblance sur le signal réel présenté dans la figure \ref{fig. exemple_signaux_reels} (c)}
\label{fig. loglik-reel}
\end{figure}

La figure \ref{fig. loglik-reel} montre  la convergence du critère de log-vraisemblance sur le signal réel présenté dans la figure \ref{fig. exemple_signaux_reels} (c). On peut remarquer pour cet exemple que, la log-vraisemblance cro\^it très rapidement jusqu'à 35 itérations puis augmente très peu jusqu'à la convergence. Pour tous les signaux traités, nous avons pu constater que l'algorithme EM convergeait en un nombre d'itérations autour de 90.


Le tableau \ref{table. resultats BIC en fonction de p} montre les pourcentages de choix  du degré $p$ des polynômes, sélectionnés par le critère BIC avec $K=5$ sous-modèles de régression qui correspondent aux phases impliquées dans une man{\oe}uvre d'aiguillage. Les valeurs $p=3$ et $p=4$ sont majoritairement sélectionnées.

 \begin{table}[!h]
\centering
\begin{tabular}{l c c c c c c }
\hline $p$ &1&2&3&4&5&6\\ \hline
\% & $0$ & $9.5238$ &  $31.4286$ &  $\boldsymbol{38.0952}$ &  $14.2857$ &   $6.6667$\\
\hline
\end{tabular}
\caption{Pourcentage de choix de $p$ avec $K=5$ obtenu sur 120 signaux réels}
\label{table. resultats BIC en fonction de p}
\end{table}

\vspace*{-0.7 cm}
\section{ Conclusion}

Une méthode de régression non linéaire a été proposée dans cet
article. Cette méthode peut être vue comme une solution alternative
au problème des moindres-carrés pour la régression non linéaire.
Elle s'appuie sur un modèle de régression simple intégrant un
processus latent qui permet d'activer successivement, et manière
souple, des sous-modèles de régression polynomiaux. Pour estimer les
paramètres du modèle proposé, un algorithme de type EM adapté au
contexte de processus latent a été proposé. La méthode proposée se distingue par le fait d'\^etre adaptée aux transitions à la fois souples et brusques gr\^ace à la modélisation particulière utilisée pour le processus latent. Dans l'étude expérimentale menée, les modèles alternatifs considérés sont un modèle de régression par morceaux et un modèle markovien de régression. Les résultats obtenus
sur des données simulées et sur des données réelles issues d'une
application du domaine ferroviaire confortent l'intérêt de
notre démarche.

\vspace*{-.4 cm}
\section*{Remerciements}

Les auteurs remercient vivement les relecteurs pour leurs remarques pertinentes, en particulier sur les intervalles de confiance en régression non linéaire.

\medskip

\bibliographystyle{rnti}

\bibliography{chamroukhi_same_govaert_aknin_rnti.bib}

\Fr

\end{document}